\documentclass[a4paper,12pt]{amsart}

\usepackage{rotating,pdflscape,color,amssymb,hyperref,subfigure,psfrag,amsmath,eufrak,bbm,epsfig}
\usepackage{a4wide,amscd}
\usepackage{tikz} 
 \usepackage[usenames,dvipsnames]{pstricks}
 \usepackage{pst-grad} 
 \usepackage{pst-plot} 

\usepackage[small]{caption}
\author[Florent~Benaych-Georges]{Florent Benaych-Georges}\address{Florent Benaych-Georges: MAP 5, UMR CNRS 8145 - Universit\'e Paris Descartes, 45 rue des Saints-P\`eres 75270 Paris cedex~6, France and CMAP, \'Ecole Polytechnique, route de Saclay, 91128 Palaiseau Cedex, France.} \email{florent.benaych-georges@parisdescartes.fr}
\author[Thierry~Cabanal-Duvillard]{Thierry Cabanal-Duvillard}
\address{Thierry Cabanal-Duvillard: MAP 5, UMR CNRS 8145 - Universit\'e Paris Descartes, 45 rue des Saints-P\`eres 75270 Paris cedex~6, France} \email{Thierry.Cabanal-Duvillard@mi.parisdescartes.fr}
\date{\today}
\thanks{This work was partially supported by the \emph{Agence Nationale de la Recherche} grant ANR-08-BLAN-0311-03 and  partly accomplished during the first named author's stay at   New York University Abu Dhabi, Abu Dhabi (U.A.E.).}

\newcommand{\bet}{\beta}
\newcommand{\Gspi}{G_{/\pi}}
\newcommand{\tGa}{\tilde{\Ga}}

\newcommand{\tm}{\tilde{m}}

\newcommand{\hk}{\hat{k}}

\newcommand{\Cc}{\mc{C}}
\newcommand{\Ga}{\Gamma} 
\newcommand{\ga}{\gamma}

\newcommand{\ld}{\ldots}
\newcommand{\beg}{\begin}
\newcommand{\en}{\end}

\newcommand{\cdd}{\cdots\cdots}

\newcommand{\trm}{\textrm}

\newcommand{\bgt}{\begin{itemize}}
\newcommand{\ent}{\end{itemize}}
\newcommand{\ite}{\item}

\newcommand{\op}{\operatorname}
\newcommand{\eqre}{\eqref}
\newcommand{\re}{\ref}
\newcommand{\la}{\label}

\newcommand{\si}{\sigma}

\newcommand{\Var}{\operatorname{Var}}

\newcommand{\lan}{\langle}
\newcommand{\ran}{\rangle}

\newcommand{\diag}{\operatorname{diag}}

\newcommand{\Part}{\operatorname{Part}}
\newcommand{\NC}{\operatorname{NC}}

\newcommand{\Pro}{\mathbb{P}}

\newcommand{\Tr}{\operatorname{Tr}}

\newcommand{\ninf}{\underset{n\to\infty}{\longrightarrow}}
\newcommand{\Ninf}{\underset{N\to\infty}{\longrightarrow}}

\newcommand{\E}{\mathbb{E}}

\newcommand{\R}{\mathbb{R}}
\newcommand{\C}{\mathbb{C}}

\newcommand{\K}{\mathbb{K}}

\newcommand{\ud}{\mathrm{d}}

\newcommand{\pro}{probability }

\newcommand{\f}{\frac}
\newcommand{\ff}{\frac{1}}
\newcommand{\lf}{\left}
\newcommand{\ri}{\right}

\newcommand{\st}{such that }

\newcommand{\lam}{\lambda}
\newcommand{\Lam}{\Lambda}
\newcommand{\ti}{\times}

\newcommand{\ste}{\, ;\, }
\newcommand{\mc}{\mathcal }

\newcommand{\eps}{\varepsilon}

\newcommand{\bxp}{\boxplus}

\newcommand{\A}{\mc{A}}

\newcommand{\NN}{\mc{N}}

\newcommand{\bck}{\backslash}
\newcommand{\al}{\alpha}
\newcommand{\tta}{\theta}

\newcommand{\eqlaw}{\stackrel{\textrm{law}}{=}}

\newcommand{\ovl}{\overline}

\newcommand{\bbm}{\begin{bmatrix}}
\newcommand{\ebm}{\end{bmatrix}}
\newcommand{\bes}{\begin{equation*}}
\newcommand{\ees}{\end{equation*}}
\newcommand{\be}{\begin{equation}}
\newcommand{\ee}{\end{equation}}
\newcommand{\beqy}{\begin{eqnarray}}
\newcommand{\eeqy}{\end{eqnarray}}
\newcommand{\beq}{\begin{eqnarray*}}
\newcommand{\eeq}{\end{eqnarray*}}
\newcommand{\one}{\mathbbm{1}}
\newcommand{\lto}{\longrightarrow}

\newcommand{\cf}{\emph{cf. }}
\newcommand{\ie}{\emph{i.e. }}
\newcommand{\eg}{\emph{e.g. }}
\newcommand{\bpm}{\begin{pmatrix}}
\newcommand{\epm}{\end{pmatrix}}

\newcommand{\cd}{\cdots}
\newcommand{\Lvy}{L\'evy }
\newcommand{\PUNmu}{\mathbb{P}_{U_N}^{(\mu)}}
\newcommand{\PtUNmu}{\mathbb{P}_{\tilde{U}_N}^{(\mu)}}

\newcommand{\enc}{\en{cases}}
\newcommand{\ensk}{\{1, \ld, k\}}
\newcommand{\thin}{\op{thin}}

\newcommand{\nc}{\op{nc}}
\newcommand{\simpi}{\stackrel{\pi}{\sim}}
\newcommand{\nsimpi}{\stackrel{\pi}{\nsim}}

\newcommand{\tGt}{\tilde{\Gamma}_t}
\newcommand{\tmu}{\tilde{\mu}}

\newtheorem{Th}{Theorem}[section]

\newtheorem{propo}[Th]{Proposition}

\newtheorem{lem}[Th]{Lemma}

\newtheorem{rmk}[Th]{Remark}

\newtheorem{propdef}[Th]{Proposition-Definition}
\newenvironment{pr}{\noindent {\bf Proof. }}{\hfill $\square$\\}
\long\def\symbolfootnote[#1]#2{\begingroup
\def\thefootnote{\fnsymbol{footnote}}\footnote[#1]{#2}\endgroup}

\newcommand{\La}{\Lambda}
\newcommand{\Card}[1]{\left\vert #1\right\vert}

\title[Marchenko-Pastur and Bercovici-Pata generalized]{Marchenko-Pastur theorem and Bercovici-Pata bijections for heavy-tailed or localized vectors}
\begin{document}
\maketitle

 \beg{abstract}The celebrated Marchenko-Pastur theorem gives the asymptotic spectral distribution of sums of random, independent,  rank-one projections. Its main hypothesis is that these projections are more or less uniformly distributed on the first grassmannian, which implies for example that the corresponding vectors are \emph{delocalized}, \ie are essentially supported by the whole canonical basis. In this paper, we propose a way to drop this delocalization assumption and we  generalize  this theorem to a quite general framework, including random projections whose corresponding vectors are \emph{localized}, \ie with some components much larger than the other ones. The first of our two main examples is given by heavy tailed random vectors (as in the model introduced by Ben Arous and Guionnet  in \cite{BAGheavytails} or as  in the model introduced by Zakharevich in \cite{Zakharevich} where the moments grow very fast as the dimension grows). Our second main example, related to the continuum between the classical and free convolutions introduced in \cite{L-BEN},  is given by vectors which are distributed as the Brownian motion on the unit sphere, with localized initial law. Our framework is in fact general enough to get new correspondences between classical infinitely divisible laws and some limit spectral distributions of random matrices, generalizing the so-called \emph{Bercovici-Pata bijection}. 
 \en{abstract}
 

\section{Introduction}

In 1967, Marchenko and Pastur introduced a successful matrix model inspired
by the elementary fact that each Hermitian matrix is the sum of orthogonal rank one
homotheties. Substituting orthogonality with independence, they considered in their seminal
paper \cite{MarchenkoPastur1967}
 the $N\times N$ random
matrix defined by
\begin{equation}\label{47117h}
  \frac{1}{N}
  \sum_{i=1}^pX^{i}\cdot U_N^{i}(U_N^i)^{*},
\end{equation}
where $(X^{i})_{i\geq 1}$ is an i.i.d.  sequence of real valued random
variables and $(U_N^{i})_{i\geq 1}$ is an i.i.d. sequence of $N$-dimensional
column vectors, whose conjugate transpose are denoted $(U^i_N)^*$, independent of $(X_i)_{i\geq 1}$. As a main result, they
proved that the empirical spectral measure of this matrix converges to a 
limit with an explicit characterization under the following assumptions:
\begin{enumerate}
\item[(a)] $N,p$ tend to
infinity in such a way that $p/N{\longrightarrow} {\lambda} >0$;
\item[(b)] the first four
joint moments of the entries of $U_N^{1}$ are not too far, roughly speaking, from the ones of the entries of a standard
Gaussian vector.
\end{enumerate}
In the special case where all $X^{i}$'s are equal
to one, the matrix (\ref{47117h}) reduces to a so-called {\it empirical
  covariance matrix}, and its limit spectral distribution is none other than the well-known \emph{Marchenko-Pastur
  distribution with parameter ${\lambda}$}.
\smallskip

It has to be noticed that even in the general case, the limit spectral distribution does not depend
  on the
particular choice of the $U_N^{i}$'s, granted they satisfy the above
hypothesis.  For example, one can choose $U_N^{1}$ to have
uniform distribution on the sphere of ${\mathbb{R}}^N$ or ${\mathbb{C}}^N$
with radius $\sqrt{N}$, or to be a standard gaussian: such a vector is said
to be \emph{delocalized}, which means that with large probability,
$\|U_N^{1}\|_\infty/\|U_N^{1}\|_2$ is small; more specifically:
$$
\frac{\|U_N^{1}\|_\infty}{\|U_N^{1}\|_2}\approx\left(\frac{\log(N)}{N}\right)^{1/2}.
$$
After the initial paper of Marchenko and Pastur, a long list of further-reaching results about limit spectral
distribution of empirical covariance matrices have been obtained, by Yin
and Krishnaiah \cite{YinKrishnaiah}, G\"otze and Tikhomirov
\cite{Gotze2004,GotzeStein}, Aubrun \cite{aubrun}, Pajor and Pastur
\cite{PajorPastur}, Adamczak \cite{adamczak}. 
All of them are devoted to
the empirical covariance matrix of more or less delocalized vectors, with
limit spectral distribution being the Marchenko-Pastur distribution (except
in the case treated in \cite{YinKrishnaiah}, but there the vector are still
very delocalized). 
\smallskip

In this paper,  our goal is to drop the delocalization assumption and to be
able to deal with \emph{localized} $U^i_N$'s, i.e. with some entries much
larger than the other ones.  For example, the applications of our main theorem   include the case where the entries of $U^i_N$ have heavy tails, but also in some other examples of localized vectors, such as the one where the law of $U_N^i$ results from a Brownian motion with localized initial condition.
\smallskip

This approach is based on our preceding works
\cite{FBDAOPinfdiv,cab-duv-BP}  on the Bercovici-Pata bijection (see also \cite{FBPTRFinfdiv,perez-abreu-Sakuma,flothierry1}). 
This bijection, that we denote by $\Lambda$, is a correspondence between the probability measures on the real line that are   infinitely divisible with respect to the classical convolution $*$ and the ones which are infinitely divisible with respect to the free convolution $\boxplus$. 
 In \cite{FBDAOPinfdiv,cab-duv-BP}, we constructed a set of matrix ensembles  which produces a quite natural
interpretation of $\Lam$. This construction is easy to describe for
compound Poisson laws and makes the connection with Marchenko-Pastur's
model quite clear. Let $(X^i)_{i\geq1}$ be still an i.i.d. sequence of real
valued random variables, $(U_N^i)_{i\geq 1}$ i.i.d.  column vectors
uniformly distributed on the sphere of $\mathbb{C}^N$ with radius
$\sqrt{N}$, and $(P(\lambda),\lambda>0)$ a standard Poisson process,
$(X^i)_{i\geq1}$, $(U_N^i)_{i\geq 1}$ and $(P(\lambda),\lambda>0)$ being
independent. For each $N\geq 1$, we defined  the random matrix
\begin{equation}\label{89111}
  \frac{1}{N}\sum_{i=1}^{P(N\lambda)}X^i\cdot U^i_N(U^i_N)^*
\end{equation}
and we proved that its empirical spectral law converges to $\La(\mu)$
when $N$ goes to infinity, where  $\mu$ is the  compound
  Poisson law of
$$
\sum_{i=1}^{P(\lambda)}X^i,
$$
  The link with Marchenko-Pastur's model of
\eqref{47117h} is now obvious and it is easy to verify that the empirical
spectral laws of  \eqref{47117h} and \eqref{89111} have same
limit if $p\sim N\lambda$. Hence, our previous works
\cite{FBDAOPinfdiv,cab-duv-BP} could be viewed as another insight on
Marchenko-Pastur's results, partly more restricted, since we only
considered uniformly distributed random vectors $U^i_N$, partly more
general, since our construction extended to all infinitely divisible laws.
In fact, the main advantage of our matricial model, which has also been studied  in \cite{perez-abreu-Sakuma},   over Marchenko-Pastur's
one is to be infinitely divisible. It allows us to derive simpler proofs,
using appropriate tools as cumulant computations or semi-groups.
\smallskip

In the present paper, we extend our construction to a larger class of
$U_N^i$'s, while continuing to benefit of the infinitely divisible
framework. Roughly speaking, we are able to prove the convergence of the
empirical spectral law if we suppose only that the entries of $U_N^i$ are
exchangeable and have a moment of order $k$ growing as $N\to\infty$ at most in $N^{\frac{k}{2}-1}$, for any fixed $k$ (\emph{cf.} Theorem \ref{mainTh5711}). 

Then, approximations  allow to extend the result to  $U_N^i$ with      \emph{heavy tailed}  entries (see Theorem \re{HT5412}).  When the $X^i$'s are constant, we recover a result by obtained by Belinschi, Dembo and Guionnet in \cite{BDGheavytails}.   Our result is   more general from a certain point of view (we allow some random $X^i$'s),  but less   explicit since we characterize the limit spectral distribution as the weak limit of a sequence of probability distributions with calculable moments and not with a functional equation, as in \cite{BDGheavytails}. We also state (Theorem \re{Zakha_MP}) a ``covariance matrices version'' of Zakharevich's generalization of Wigner's theorem of \cite{Zakharevich},  which is a direct consequence of our key result Theorem \re{mainTh5711}.
\smallskip

We also devote a particular interest to the special case where the column
vectors $ U_{N}^i$'s are copies of
\begin{equation}\label{2191121h47}
  U_N(t)   =\sqrt{N}e^{-\frac{t}{2}}X_N+ \sqrt{1-e^{-t}} G_N,
\end{equation}
where $X_N$ is uniformly distributed on the canonical basis $(1, 0,
{\ldots}, 0)^T$,{\ldots}, $(0, {\ldots}, 0,1)^T$, independent of $G_N:=(z_1, {\ldots} ,z_N)^T$, with $z_1, z_2,z_3,
{\ldots}$ independent standard   Gaussian variables.  Let us
emphasized that $U_N(t)$ is a quite typicaly \emph{localized} vector since
it has exactly one entry which is much bigger than the $N-1$
others. Precisely, we have:
$$
\frac{\|U_N(t)\|_\infty}{\|U_N(t)\|_2}\approx e^{-t/2}.
$$
This is the reason why the $N\to\infty$ limits differ from the ones of the
classical matrix models. For example, the classical Marchenko-Pastur
Theorem about empirical covariance matrices is not true anymore for such
vectors.  More specifically, we have:

\begin{Th}
  Let $U_N^1, U_N^2, {\ldots}$ be i.i.d. copies of the vector $U_N(t)\in
  {\mathbb{C}}^{N}$ defined in \eqref{2191121h47}. Let
  \begin{equation}\label{Mfig_hang_on}
    M=\frac{1}{N}\sum_{i=1}^pU_N^i( U_N^i)^*
  \end{equation}
  be the (dilated) \emph{empirical covariance matrix} of the sample $U_N^1,
  {\ldots}, U_N^p$. As $N,p\to\infty$ with $p/N\to{\lambda}>0$, the
  empirical spectral law of $M$ converges to a limit law
  $P_{{\lambda},t}$ with unbounded support which is characterized by its
  moments, given by the formula
  \begin{equation}\label{21911.22h35}
    \int x^kP_{{\lambda},t}({\mathrm{d}} x)=\sum_{\pi\in{\operatorname{Part}}(k)}e^{-\kappa(\pi)t} {\lambda}^{|\pi|}
  \end{equation}
  with ${\operatorname{Part}}(k)$ the set of partitions of $\{1,\ldots,k\}$ and $\kappa$ defined by:

  $\bullet$ $\kappa(\pi)=\kappa(\pi_1)+\cdots+\kappa(\pi_q)$ if
  $\pi_1,\ldots,\pi_q$ are the connected components of $\pi$;

  $\bullet$ if $\pi$ is connected, $\kappa(\pi)$ is the number of times one
  changes from one block to another when running through $\pi$ in a cyclic
  way.
\end{Th}

Note that the law $P_{{\lambda},t}$ is the Poisson law with parameter
${\lambda}$ if $t=0$ and tends to the free Poisson law (also called
Marchenko-Pastur law) as $t$ tends to $+\infty$: the weight
$e^{-\kappa(\pi)t}$ penalizes {\it crossings} in partitions (indeed,
$\kappa(\pi)=0$ {if and only if } $\pi$ is non crossing). An illustration
is given in Figure \ref{fig_MP_del} below. In fact, Formula
\eqref{21911.22h35} can be generalized into a more general
moments-cumulants formula which provides a new continuous interpolation
between classicaly and freely infinitely divisible laws. This interpolation is related to the notion of $t$-freeness developed by the first named author and L\'evy in \cite{L-BEN} and is based on a progressive penalization of the crossings in the moments-cumulants formula  \eqre{741216h48}.

\begin{figure}[h!]
  \centering \subfigure[Case where $t=0.01$]{
    \includegraphics[height=1.7in, width=3.05in]{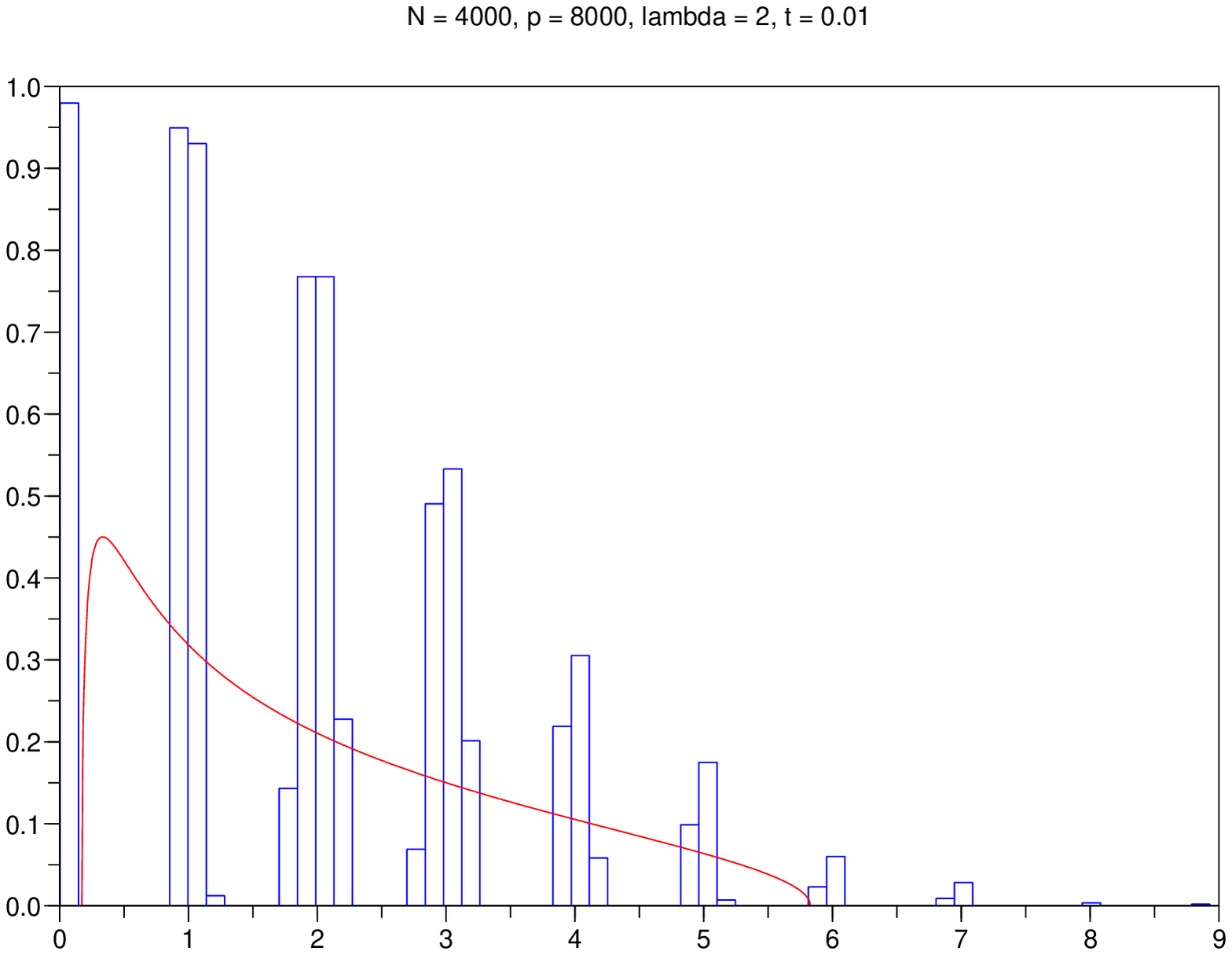}
    \label{MPdeloc_t0.01}} \subfigure[Case where $t=0.1$]{
    \includegraphics[height=1.7in, width=3.05in]{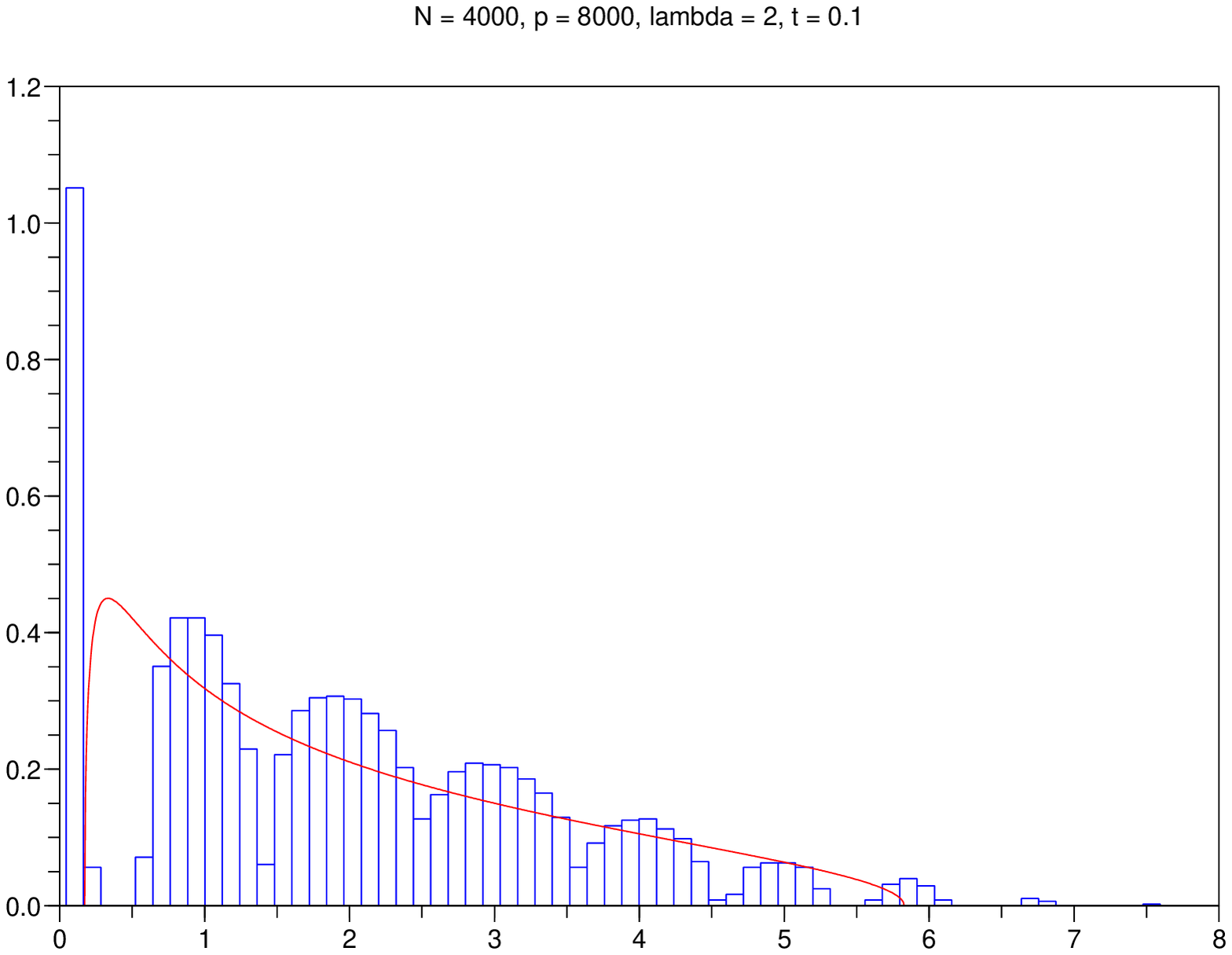}
    \label{MPdeloc_t0.1}} \subfigure[Case where $t=1$]{
    \includegraphics[height=1.7in, width=3.05in]{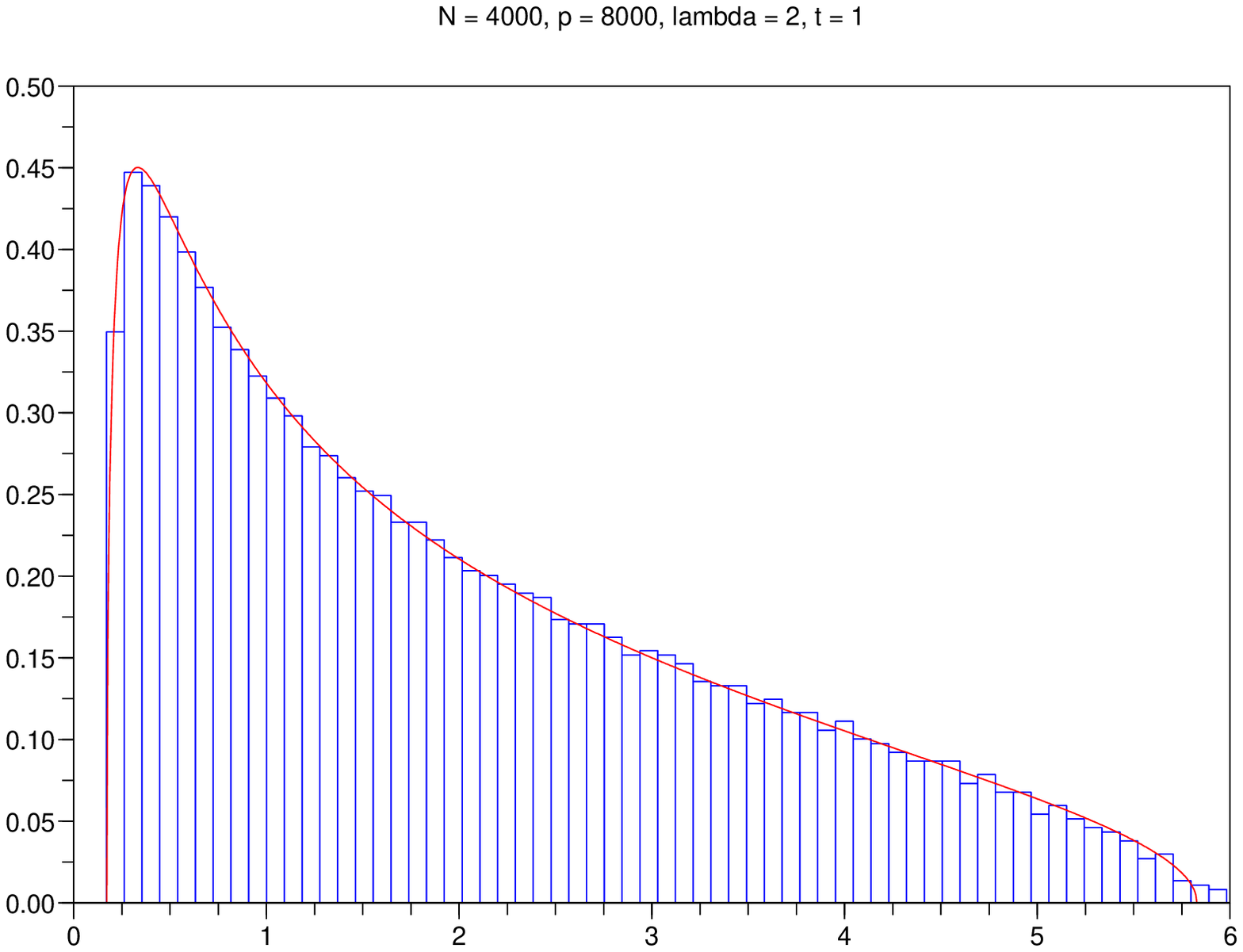}
    \label{MPdeloc_t1}}
  \caption{The empirical spectral distribution of the matrix $M$ of
    \eqref{Mfig_hang_on} for $N=4.10^3$, $p=8.10^3$ at several values of
    $t$. We see that this distribution, which is approximately the law
    $P_{2,t}$, is close to the {\bf Poisson law} with parameter $2$ for $t$ close
    to zero, and that it converges to the {\bf Marchenko-Pastur law} with
    parameter $2$ (whose density is plotted by a smooth continuous line) as $t$ grows.}
  \label{fig_MP_del}
\end{figure}


The paper is organized as follows. Next section describes the main results;
Section 3 provides some applications; the
following sections are devoted to the proofs; the last section is an
appendix where we recall some facts on the infinitely divisible laws, the
Bercovici-Pata bijection and (hyper)graphs, for the easyness of the reader.

As this preprint was published, Victor P\'erez-Abreu
informed us that he is working on a close subject with J. Armando
Dom\'\i nguez-Molina  and Alfonso Rocha-Arteaga in the forthcoming preprint  \cite{vpajadara}.

\medskip

\paragraph{{\bf Acknowledgements.}} The authors thank James
Mingo for a useful explanation during the program ``Bialgebras in free
probability'' in the Erwin Schr\"odinger Institute in Wien during spring
2011. They also thank Thierry {L\'evy} and Camille Male for some discussions on a preliminary
version of the paper.

\section{The main results} 

\subsection{A general family of matrix ensembles}

  The basic facts on the \emph{$*$-infinitely divisible} laws
are recalled in the appendix.  For $\mu$ such a law, its \emph{\Lvy
  exponent} is the function $\Psi_\mu$ \st the Fourier transform of $\mu$
is $e^{\Psi_\mu(\cdot)}$. 
\smallskip

\subsubsection{Case of  compound Poisson laws}\la{COCPL5312}

Such a law $\mu$ is
a one of $$\sum_{i=1}^{P(\lam)}X_i,$$ where $P(\lam)$ is a random variable
with Poisson law with expectation $\lam$ and the $X_i$'s are i.i.d. random
variables, independent of $P(\lam)$. In this case, if $\nu$ denotes the law
of the $X_i$'s, the \Lvy exponent of $\mu$
is 
$
\Psi_\mu(\xi)=\lam\int(e^{it\xi}-1)\ud\nu(t).
$ For $N\ge 1$, let $(U_N^{i})_{i\ge 1}$ be a sequence of i.i.d. copies of  a random column vector $U_N\in \R^{N\ti 1}$ or $\C^{N\ti 1}$. Then we define   $\PUNmu$ to be  the law
of $$\ff{N}\sum_{i=1}^{P(N\lam)}X_i\cdot U_N^i(U_N^i)^{*},$$ where $P(N\lam)$ is a random variable with Poisson law with expectation $N\lam$, the $X_i$'s are i.i.d. random variables with law $\nu$ and the $U_N^i$'s are i.i.d. copies of 
$U_N$ (whose conjugate transpose are denoted by $(U_N^i)^{*}$), all being independent. Let us notice that $\PUNmu$ is still a compound Poisson law and that its L\'evy exponent is given by:
$$
A\mapsto N\E\left[\Psi_\mu(U_N^*AU_N/N)\right]
$$
for any $N\times N$ Hermitian matrix $A$.

\subsubsection{General case} Here, we shall extend the previous construction for  $\mu$   a general infinitely divisible law.
In the following theorem and in the rest of the paper, the spaces of \pro measures are endowed with the weak topology. $\K$ denotes either $\R$ or $\C$.

\begin{Th}\la{defPUNmu}
Let $\mu$ be an infinitely divisible law on $\R$, let us fix $N\ge 1$ and let $U_N\in \K^{N\ti 1}$ be
  a  random column vector such that $\E\left[\Vert U_N\Vert_2^4\right]<+\infty$. Let $(U_N^{i})_{i\ge 1}$ be a sequence of i.i.d. copies of $U_N$ and, for each $n\ge 1$, let  $X_n^{1},  X_n^{2}, \ld, X_n^{Nn} $ be i.i.d. $\mu^{*\ff{n}}$-distributed random variables. Then the sequence of $N\ti N$ Hermitian matrices   \be\la{47117h1prime}\ff{N}\sum_{i=1}^{N\ti  n}X_n^{i}\cdot U_N^{i}(U_N^{i})^*\ee converges in distribution, as $n\lto\infty$ ($N$ being fixed), to a \pro measure $\PUNmu$ on the space of $N\ti N$ Hermitian matrices, whose Fourier transform   is    given by \be\la{Fourier_loi6711}\int e^{i\Tr(AM)}\ud\PUNmu(M)=\exp\lf\{N\E\left[\Psi_\mu(U_N^*AU_N/N)\right]\ri\}\ee for any $N\ti N$ Hermitian matrix $A$.
 \en{Th}

  The following proposition extends the theorem to a quite more general framework, where the $X_n^i$'s are i.i.d. but not necessarily distributed according to $\mu^{*\ff{n}}$ and only satisfy the limit theorem $$\trm{law of }(X_n^1+\cd\cd+X_n^{k_n})\;\ninf \;\mu$$ for a sequence $k_n\to\infty$.

\beg{propo}\la{129111} a)  If the law of  $U_N$ is compactly supported, then  to any limit theorem $$
\underbrace{\nu_n*\cdd*\nu_n}_{k_n\trm{ times}}\ninf \mu,
$$ there corresponds a limit theorem  in the space of $N\ti N$ Hermitian matrices   \be\la{47117h1}\ff{N}\sum_{i=1}^{N\ti k_n}X_n^{i}\cdot U_N^{i}(U_N^{i})^*\ninf \PUNmu,\ee where $X_n^{1},  X_n^{2}, \ld , X_n^{Nk_n}$ are i.i.d. $\nu_n$-distributed random variables, $(U_N^{i})_{i\ge 1}$ is a sequence of i.i.d. copies of $U_N$, independent of the $X_n^{i}$'s.

b) In the case where the law of $U_N$ is not compactly supported, \eqre{47117h1} stays true as long as one supposes that $k_n\ti \int |t|\nu_n(\ud t)$ is bounded uniformly in  $n$.
\en{propo}

\beg{rmk}{\rm  Such a construction has been generalized to the
more general setting of Hopf algebras by Sch\"urmann, Skeide and Volkwardt in \cite{2007arXiv0712.3504S}.}\en{rmk}

\subsection{Convergence of the empirical spectral  law}
The \emph{empirical spectral law} of a matrix is the uniform \pro measure on its eigenvalues (see Equation \eqre{710129h45} in the appendix).
When $U_N/\sqrt{N}$ is uniformly distributed on the unit sphere, we proved
that the empirical spectral law associated to $\PUNmu$ converges when the
size $N$ tends to infinity (\emph{cf.} \cite{FBDAOPinfdiv,cab-duv-BP}). In
order to obtain other convergences, we shall first make the following
assumptions on the random vector $U_N\in \K^{N\ti 1}$.  We denote the entries of $U_N$ by $U_N(1), \ld, U_N(N)$. Roughly speaking, these assumptions  mean that the $U_N(i)$'s are exchangeable and that the moment of order $k$ of $U_N(1)$ grows at most in $N^{\f{k}{2}-1}$, for all $k$. These assumptions will be weakened in the next section to consider heavy tailed variables.

\beg{hyp}\la{hyp_NfixeUn} a) For each    $N\ge 1$, the entries of $U_N$ are exchangeable and  have moments of all orders.

b) As $N$ goes to infinity, we have: \bgt
\ite 
  for each $p\ge 1$, for each $i_1, \ld, i_{2p} \ge 1$,   \be\la{571117h30_mjj}\E[U_N(i_1)\cd U_N(i_p)\ovl{U_N(i_{p+1})}\cd\ovl{U_N(i_{2p})}]=O(N^{p-|\{i_1, \ld, i_{2p}\}|}),\ee
  \ite
  for each $k\ge 1$, for all positive integers $ n_1,  \ld, n_k$,  there exists $\Gamma(n_1, \ld, n_k)$ finite \st    
\be\la{571117h30}
\E\left[\f{|U_N(1)|^{2n_1}}{N^{n_1-1}}\cd\cd\f{|U_N(k)|^{2n_k}}{N^{n_k-1}}\right]\Ninf \Gamma(n_1,
\ld, n_k).
\ee Moreover,  there is a constant $C$ \st for all $k\ge 1$, for all  $n_1, \ld, n_k$, \be\la{571117h30ccbain}\Gamma(n_1,
\ld, n_k)\le C^{n_1+\cd+n_k}.\ee
\ent
\en{hyp}

\beg{rmk}\la{tildeU_N}{\rm Let define the random column vector $\tilde{U}_N=\diag(\eps_1, \ld, \eps_{N})U_N$, with $\eps_1, \ld, \eps_{N}$ some i.i.d. variables   uniformly distributed  on $\{z\in \K\ste |z|=1\}$, independent of $U_N$. Note that if $U_N$ satisfies  Hypothesis 
\re{hyp_NfixeUn},   so does $\tilde{U}_N$, with the same function $\Ga$. Moreover,   the expectation of \eqre{571117h30_mjj} is null for $\tilde{U}_N$   as soon as the multisets $\{i_1, \ld, i_k\}$ and $\{i_{k+1}, \ld, i_{2k}\}$ are not equal (in the complex case) and as soon as an element in the multiset $\{i_1, \ld, i_{2k}\}$ appear an odd number of times (in the real case).}\en{rmk}
  

The   following theorem 
is the key result of the paper. For its second part, we endow the set of   functions on multisets of integers (it is the set that $\Ga$ belongs to) with the product topology. 
\beg{Th}\la{mainTh5711}
We suppose that Hypothesis  \re{hyp_NfixeUn} holds.

a) For any $*$-infinitely divisible distribution $\mu$, the empirical spectral distribution of a $\PUNmu$-distributed random matrix converges almost surely, as $N\lto\infty$, to a deterministic \pro measure $\Lam_\Gamma(\mu)$, which depends only on $\mu$ and on the function $\Gamma$ of Equation \eqre{571117h30}. 

b) The \pro measure $\Lam_\Gamma(\mu)$ depends continuously on the pair $(\mu, \Ga)$. 

c) The moments of $\Lam_\Gamma(\mu)$
moments can be computed when $\mu$ has moments to all orders via the following formula,
\be\la{eq:0519-2ter}
 \qquad\qquad\qquad\qquad \int x^k \Lam_\Gamma(\mu)(\ud x)\;=\sum_{\pi\in \Part(k)}f_\Ga(\pi)\prod_{J\in \pi} c_{|J|}(\mu),\qquad\qquad (k\ge 1)
\ee
where the non negative numbers $f_\Ga(\pi)$, which factorize along the connected components of $\pi$,  are given at Lemma \re{pte:0213-3gusgus} and the numbers $c_{|J|}(\mu)$ are the cumulants of $\pi$ (whose definition is recalled in  the appendix).  

d) If the \Lvy measure of $\mu$ has compact support, then $\La_\Ga(\mu)$ admits exponential moments of all orders.
\end{Th}



The following proposition 
allows to assert that many  limit laws obtained in Theorem \re{mainTh5711} have unbounded support. Recall that a law is said to be {\it non degenerate} if it is not a Dirac mass.
\beg{propo}\la{unbounded_support_121011}If the function $\Ga$ of Hypothesis  \re{hyp_NfixeUn} is \st $\inf \Ga(n)^{1/n}>0$, then for any non degenerate 
  $*$-infinitely divisible distribution $\mu$ whose classical cumulants are all non negative, the law $\Lam_\Ga(\mu)$ has unbounded support.  
\en{propo}

\section{Applications  and examples}

\subsection{Heavy-tailed Marchenko-Pastur theorem}
In this section, we   use Theorem \re{mainTh5711} to extend the theorem of Marchenko and Pastur described in the introduction to vectors with heavy-tailed entries. This  also extends Theorem 1.10 of Belinschi, Guionnet and Dembo in \cite{BDGheavytails} (their theorem   corresponds to the case $X_i=1$), but we do not provide the explicit characterization of the limit that they propose.  Also, our result is valid in the complex as in the real case, whereas the one of \cite{BDGheavytails} is only stated in the real case (but we do not know whether this is an essential restriction of the approach of \cite{BDGheavytails}).

Let us fix $\al\in (0,2)$ and let $(Y_{ij})_{i,j\ge 1}$ be an infinite  array of i.i.d.  $\K$-valued random variables  \st the function \be\la{queue_distrib_6412}L(t):= t^\al \mathbb{P}(|Y_{11}|\ge t)\ee has slow variations as $t\to\infty$. 
For each $j$, define the column vector  $V^j:=(Y_{ij})_{i=1}^N\in\K^{N\ti 1}$ ($V^j$ depends implicitly on $N$).
Let us define $a_N:=\inf\{t\ste t^{-\al}L(t)\le 1/N\}$.
\beg{Th}\la{HT5412}For any set $(X_j)_{j\ge 1}$ of i.i.d. real random variables (with any law, but that does not depend on $N$)  independent of the $Y_{ij}$'s and any fixed $\lam>0$, the empirical spectral law of the random matrix $$M_{N,p}:=\ff{a_N^2}\sum_{j=1}^p X_j\cdot V^j(V^j)^*$$ converges almost surely, as $N,p\to\infty$ with $p/N\to \lam$, to a deterministic \pro measure which depends only on $\al$, on $\lam$ and on the law of the $X_i$'s.  This limit  law depends continuously of these three parameters.
\en{Th}

\subsection{Covariance matrices with exploding moments}
The following theorem is a direct consequence of Theorem \re{mainTh5711} and Proposition \re{unbounded_support_121011}. It is the ``covariance matrices version'' of Zakharevich's generalization of Wigner's theorem to   matrices with exploding moments  (see \cite{Zakharevich} and also the recent work \cite{CamilleTrafics} by Male).  

\beg{Th}\la{Zakha_MP}Let $M_{N,p}=[X_{ij}]$ be a complex $N\ti p$ random matrix with i.i.d. centered entries whose distribution might depend on $N$ and $p$. We suppose that there is a sequence $\mathbf{c}=(c_k)_{k\ge 2}$  \st $c_k^{1/k}$ is bounded and \st for each fixed $k\ge 2$, $$\qquad\qquad\qquad\qquad\f{\E[|X_{11}|^k]}{N^{\f{k}{2}-1}}\lto  c_k\qquad\qquad\trm{ as   $N,p\to\infty$ with $p/N\to\lam>0$}.$$Then the empirical spectral law of $$\ff{N}M_{N,p}M_{N,p}^*$$ converges, as   $N,p\to\infty$ with $p/N\to\lam>0$, to a \pro measure $\mu_{\lam,\mathbf{c}}$ which depends continuously on the pair $(\lam, \mathbf{c})$. If $c_k=0$ for all $k\ge 3$, $\mu_{\lam,\mathbf{c}}$ is  Marchenko-Pastur distribution with parameter $\lam$ dilated by a coefficient  $\sqrt{c_2}$. Otherwise,   $\mu_{\lam,\mathbf{c}}$ has unbounded support but admits exponential moments of all orders.  
\en{Th}

\subsection{Non centered Gaussian vectors, Brownian motion on the unit sphere and a continuum of Bercovici-Pata's bijections}
 
In this section, we give  examples of vectors $U_N$ which satisfy Hypothesis \re{hyp_NfixeUn}  for a certain function $\Ga$. These examples will allow to construct the continuum of Bercovici-Pata bijections $(\Lam_t)_{t\ge 0}$ mentioned in the introduction. 
 
Let us fix $t\ge 0$ and  consider $U_N=U_{N,t}$ defined thanks to one of the two (or three, actually)  following definitions \beg{enumerate}\ite[(a)] {\it Gaussian and uniform cases}:   \be\la{defGauss_111011_71012}U_{N,t}:=e^{-\f{t}{2}}\sqrt{N}X_N+ \sqrt{1-e^{-t}} G_N,\ee
where $X_N$ is uniformly distributed on the canonical basis $(1, 0, \ld, 0)^T$,\ld $(0, \ld, 0,1)^T$  of $\K^{N\ti 1}$, independent of  $G_N:=(g_1, \ld ,g_N)^T$, with $g_1, g_2,g_3, \ld$ i.i.d.  variables whose distribution is either the centered  Gaussian law on $\K$  with variance $1$ or the uniform law on $\{z\in \K\ste |z|=1\}$,\\ 
\ite[(b)] {\it Brownian motion:} $U_{N,t}/\sqrt{N}$ is a {\it Brownian motion on the unit sphere of $\C^{N\ti 1}$}  taken at time $t$, whose distribution at time zero is the uniform law on the canonical basis of $\C^{N\ti 1}$. Such a process   is    a strong solution of the SDE \be\la{EDSbro}\ud U_{N,t}=(\ud K_t)U_{N,t}-\ff{2}U_{N,t}\ud t,\ee where $K=(K_t)_{t\ge 0}$ is a standard Brownian motion on the Euclidian  space of $N\ti N$ skew-Hermitian matrices, endowed with the scalar product $A\cdot B=N\Tr(A^*B)$.
 \en{enumerate}
Of course, these models make sense for $t=+\infty$ : in the first model, it means only that $U_{N,t}=G_N $ and the second model, at $t=+\infty$, can be understood as its limit in law, \ie a random vector with uniform law on the sphere with radius $\sqrt{N}$. For $t=+\infty$, all formulas below make sense (and stay true) by taking their $t\to +\infty$ limits. 

Thanks to the results of \cite{flobulletin}, it can be seen that the Gaussian model and the Brownian one have approximately the same finite-dimensional marginals, as $N\to\infty$. The following proposition, whose proof is postponed to Section \re{871119h31}, makes this analogy stronger. 

\beg{propo}\la{871119h30}
For $U_N=U_{N,t}$ as defined according to any of the above models,   Hypothesis    \re{hyp_NfixeUn}  holds for $\Ga:=\Ga_t$ given by the following formulas:\be\la{871119h29gauss}
\Ga_t(\underbrace{1, \ld,1}_{k\trm{ times}})=(1-e^{-t})^k+ke^{-t}(1-e^{-t})^{k-1}\quad\trm{ for every $k\ge 1$,}\ee\be
\la{871119h292gauss}\Ga_t(n, \underbrace{1, \ld,1}_{k\trm{ times}})=e^{-nt}(1-e^{-t})^k\quad\trm{ for every $n\ge 2$ and $k\ge 0$,}\ee\be
\la{871119h293gauss}\Ga_t(n_1, \ld,n_k)=0\quad\trm{ if there is $i\ne j$ \st $n_i\ge 2$, $n_j\ge 2$.}\ee
\en{propo}


Let us now consider the family of transforms $(\Lambda_{\Gamma_t},t\geq0)$
(denoted by $\Lambda_t$ in the introduction).

In both models, when $t=0$, the rank-one random projector $U_{N,t}U_{N,t}^*$ is a diagonal
matrix with   unique non zero entry uniformly distributed on the
diagonal and equal to $N$. In such a case, it can be readily seen that
$\mathbb{P}_{U_{N,0}}^{(\mu)}$ is the law of a diagonal matrix with
i.i.d. entries with law $\mu$. Owing to the law of large numbers,
$\Lambda_{\Gamma_0}$ is obviously the identity map. Moreover, it has been seen in \cite{FBDAOPinfdiv,cab-duv-BP} that for $t=\infty$, $\Lam_{\Ga_\infty}$ is the Bercovici-Pata bijection.


Therefore,  $\left(\Lambda_{\Gamma_t},t\in [0,+\infty]\right)$ provides a
continuum of   maps passing from the identity ($t=0$) to the
Bercovici-Pata bijection ($t=+\infty$). This continuum is related to the notion of $t$-freeness developed by the first named author and L\'evy in \cite{L-BEN}. The maps $\Lambda_{\Gamma_t}$ are made explicit (at least at the moments level) by the following proposition, whose proof, based on the explicitation of the functions $f_{\Ga_t}$, is postponed to   Section \ref{sec:kappa}. 

\begin{propo}\la{propo:kappa} 
Let $\mu$ be an infinitely divisible law with moments of all orders. Then
for each $t\geq 0$, 
\be\la{741216h48}\qquad\qquad\qquad\qquad
\int
x^k\Lambda_{\Gamma_t}(\mu)(\mathrm{d}x)\;=\sum_{\pi\in\Part(k)}e^{-\kappa(\pi)t}c_{\pi}(\mu)\qquad\qquad(k\ge 1)
\ee
with $\kappa$ defined by:
 
  $\bullet$ $\kappa(\pi)=\kappa(\pi_1)+\cdots+\kappa(\pi_q)$ if $\pi_1,\ldots,\pi_q$ are the connected components of $\pi$;

$\bullet$ if $\pi$ is connected, $\kappa(\pi)$ is the number of times one changes from one block to another when running through $\pi$ in a cyclic way.\end{propo}

\begin{figure}[h!]
\scalebox{1} 
{
\begin{pspicture}(0,-0.445)(7.6596875,2.82)
\psdots[dotsize=0.2](0.1,0.0)
\psdots[dotsize=0.2](0.9,0.0)
\psdots[dotsize=0.2](1.7,0.0)
\psdots[dotsize=0.2](2.5,0.0)
\psdots[dotsize=0.2](3.3,0.0)
\psdots[dotsize=0.2](4.1,0.0)
\psdots[dotsize=0.2](4.9,0.0)
\psdots[dotsize=0.2](5.7,0.0)
\psdots[dotsize=0.2](6.5,0.0)
\psdots[dotsize=0.2](7.3,0.0)
\usefont{T1}{ptm}{m}{n}
\rput(0.10375,-0.29){1}
\usefont{T1}{ptm}{m}{n}
\rput(0.9671875,-0.29){2}
\usefont{T1}{ptm}{m}{n}
\rput(1.7453125,-0.29){3}
\usefont{T1}{ptm}{m}{n}
\rput(2.5718749,-0.29){4}
\usefont{T1}{ptm}{m}{n}
\rput(3.3490624,-0.29){5}
\usefont{T1}{ptm}{m}{n}
\rput(4.160625,-0.29){6}
\usefont{T1}{ptm}{m}{n}
\rput(4.960625,-0.29){7}
\usefont{T1}{ptm}{m}{n}
\rput(5.7475,-0.29){8}
\usefont{T1}{ptm}{m}{n}
\rput(6.5603123,-0.29){9}
\usefont{T1}{ptm}{m}{n}
\rput(7.425,-0.29){10}
\psarc[linewidth=0.04](2.899889,0.0){2.8000002}{-1.5481567}{180.0}
\psarc[linewidth=0.04](6.499889,0.0){0.8}{0.0}{180.0}
\psarc[linewidth=0.04](1.7,0.0){0.8}{0.0}{180.0}
\psarc[linewidth=0.04](2.5,0.0){0.8}{0.0}{180.0}
\psarc[linewidth=0.04](5.7000003,0.04){0.8}{0.0}{180.0}
\psarc[linewidth=0.04](4.5,0.08){0.4}{0.0}{180.0}
\end{pspicture} 
}
\caption{The partition $\pi=\left\{\{1,8,10\},\{2,4\},\{3,5\},\{6,7,9\}\right\}\in\Part(10)$.}\label{fig:pi}
\end{figure}
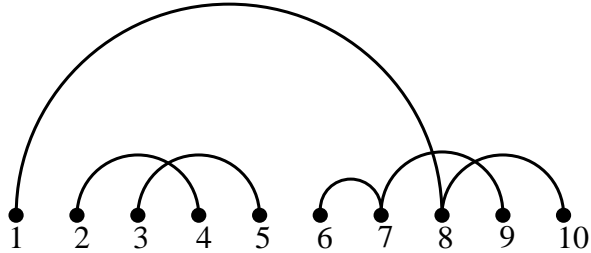

\paragraph{{\bf Example of computation of $\kappa(\pi)$.}} 
For instance, let us consider the partition
$\pi=\left\{\{1,8,10\},\{2,4\},\{3,5\},\{6,7,9\}\right\}\in\Part(10)$ (\emph{cf.}
Figure \ref{fig:pi}). Then the connected components of $\pi$ are the partitions $\pi_1$ and $\pi_2$ induced by $\pi$ on the sets $\{1,6,7,8,9,10\}$ and $\{2,3,4,5\}$. Since $\kappa(\pi_1)= 5$ and $\kappa(\pi_2)=4$, $\kappa(\pi)=9$.  

Let us now say a few words about the bijection $\Lam_{\Gamma_t}$, as a function of $t\in [0, +\infty)$.

If $t=0$, then:
$$
\int
x^k\Lambda_{\Gamma_0}(\mu)(\mathrm{d}x)=\sum_{\pi\in\Part(k)}c_{\pi}(\mu)=\int
x^k\mu(\mathrm{d}x).
$$
This corroborates the fact that  $\Lambda_{\Gamma_0}$ is the identity map.

If $t$ tends to infinity, then:
$$
\lim_{t\rightarrow+\infty}\int
x^k\Lambda_{\Gamma_t}(\mu)(\mathrm{d}x)=\sum_{\pi\in\NC(k)}c_{\pi}(\mu)=\int
x^k\Lambda({\mu})(\mathrm{d}x)
$$
where $\NC(k)$ denotes the set of non-crossing partitions of $\{1,\ldots,k\}$ (see Section \ref{sect_partitions_19911} for a precise definition). Indeed, $\kappa(\pi)$ is positive unless $\pi$ is non-crossing. Since the
continuity of $\Lambda_{\Gamma}$ w.r.t. the weak topology is uniform in
$\Gamma$, as it appears from the proof of Proposition \ref{571119h24} b),
this proves that $\Lambda_{\Gamma_t}$ tends to the Bercovici-Pata bijection $\Lambda$
when $t$ goes to infinity, as expected.


\subsection{The distribution $\PUNmu$ made more explicit} When $\mu$ is a compound Poisson law, the definition of $\PUNmu$ has been made explicit in Section \re{COCPL5312}. In this section, we explicit $\PUNmu$ for $\mu$ a    Dirac mass or a  Gaussian laws   (the column vector $U_N$ underlying the definition of $\PUNmu$ staying as general as possible). Since any
infinitely divisible law is a weak limit of convolutions of such laws and
obviously, by the Formula \eqre{Fourier_loi6711} of the Fourier transform
of $\PUNmu$, the law $\PUNmu$ depends continuously on $\mu$ and
satisfies 
$$
\PUNmu*\mathbb{P}_{U_N}^{(\nu)}=\mathbb{P}_{U_N}^{(\mu*\nu)},
$$
this gives a good idea of what a random matrix distributed according to
$\PUNmu$ looks like. Moreover, we also consider the case where $\mu$ is a
Cauchy law, where a surprising behaviour w.r.t. the convolution appears. 


First, it can easily be seen that if $\mu$ is the Dirac mass at  $\ga\in \R$, then $\PUNmu$ is the Dirac mass at 
$$
\ga \E[U_NU_N^*]=\gamma\left(\mathbb{E}\left[\vert U_N(1)\vert^2\right]-\mathbb{E}\left[ U_N(1)\overline{U_N(2)}\right]\right)+\mbox{ rank-one matrix}.
$$
Hence, due to Hypothesis \ref{hyp_NfixeUn} and Lemma \ref{lemCV6312},   
$\Lam_\Ga(\delta_\ga)=\delta_{\ga'}$ for  $\ga'=\ga\Ga(1)$.

Suppose that now $\mu$ is the standard Gaussian law $\NN(0,1)$ and that Hypotheses 
\re{hyp_NfixeUn}  
holds. Let $\tilde{U}_N$ be as in Remark \re{tildeU_N}.
Then the distribution $\PtUNmu$ only depends on $\al:=\E[|U_N(1)|^2|U_N(2)|^2]$,  $\ga:=\E\left[\vert
  U_N(1)\vert^4\right]$ and $\bet:=\dim_\R(\K)$ : when $\ga\ge \al$,  $\PUNmu$ is the law of $$\sqrt{\f{\ga-2\al}{N}}\bpm g_1&&\\ &\ddots&\\ &&g_N\epm+ \sqrt{\f{\al}{N}}\bpm g&&\\ &\ddots&\\ &&g\epm +\sqrt{\f{2\al}{N\bet}} \mathrm{GO(U)E},$$
where $g_1, \ld, g_N, g$ are standard Gaussian variables and $ \mathrm{GO(U)E}$ designs a GOE or GUE matrix (according to weither $\K=\R$ of $\C$) as defined p. 51 of \cite{agz09} (\ie a standard Gaussian vector on the space of real symmetric or Hermitian matrices endowed with the scalar product $A\cdot B:=\f{\bet}{2}\Tr AB$).
As a consequence, since when $N\to\infty$,  $\al \sim \Ga(1,1)$ and $\ga\sim N\Ga(2)$,   the spectral law of a $\PtUNmu$-distributed matrix converges to  $$\Lam_\Ga(\mu)=\NN(0, \Ga(2))\bxp (\trm{Semicircle law with variance $2\Ga(1,1)/\bet$}).$$

In the next proposition, we consider the case where $\mu$ is a  {\it Cauchy law with paramater $t$}:
$$
\mu=\Cc_t(\ud x):= \f{t\ud x}{\pi(t^2+x^2)}.
$$ 
Let $P_t$
be the associated kernel, defined by 
$$
P_t(f)(x)=\int  f(x+y)\Cc_t( \ud y).
$$

\beg{propo}\la{771113h50}
 We
suppose   that $\E\left[U_NU_N^*\right]=I_N$. For   $M_t$ a  $\mathbb{P}_{U_N}^{(\Cc_t)}$-distributed random matrix and   $f : \R\to \R$ bounded, for any $N\ti N$ Hermitian matrix
$A$, 
$$
\E\left[f(A+M_t)\right]=P_{t}(f)(A),
$$ 
real functions being applied to Hermitian matrices via the functional calculus. 
\en{propo}

  The following Corollary  follows easily, using the vector $V_N:=\tilde{U}_N/\sqrt{\Ga(1)}$  instead of $U_N$ (with $\tilde{U}_N$ as defined in Remark \re{tildeU_N}).

\beg{cor}\la{Cauchy_fixed_101011} Under  Hypothesis 
\re{hyp_NfixeUn}, the set of  Cauchy laws is invariant by $\Ga_\Lam$ : for any $t>0$,  $\Lam_\Ga(\mc{C}_t)=\mc{C}_{t'}$, for $t'=\Ga(1)t$.
\en{cor}

\section{Proof of Theorem \re{defPUNmu} and Proposition \re{129111}}\la{prdefPUNmu}
One proves the theorem and the proposition in the same time, by showing  that under the hypotheses of the theorem or of a) or b) of the remark, the Fourier transform of the  left hand term of \eqre{47117h1}, that we denote by $M_n$,  converges pointwise to the right  hand term of \eqre{Fourier_loi6711}. Indeed, for any Hermitian matrix $A$, 
$$\E[e^{i\Tr(AM_n)}]=\lf(1+\f{\E[\Psi_n(U_N^*AU_N/N)]}{k_n}\ri)^{k_n},$$
where for any $x\in \R$, $\Psi_n(x):=k_n\int (e^{itx}-1)\ud\nu_n(t)$.
The function $\Psi_n$ converges to $\Psi_\mu$ as $n\lto\infty$, uniformly
  on every compact subset of $\R$ (this follows   from  \cite[Lem. 3.1]{petrov} and from the fact that $\hat{\nu}_n^{k_n}\lto\hat{\mu}$ uniformly on every compact set). It is enough to prove the result in the case where  the law of $U_N$ is compactly supported. To conclude in the case where $\nu_n=\mu^{*\ff{k_n}}$ (resp. in the case where $\int |t|\nu_n(\ud t)=O(1/k_n)$), one needs to argue that  there is a constant $C$ \st for all $t>0$,  $\Psi_{\mu^{t}}(\lambda)\leq
  tC\left(\vert\lambda\vert+\lambda^2\right)$ (resp.   that $|\Psi_n(x)|\le k_n \int |t|\nu_n(\ud t)$).

\section{Preliminaries for the proof of Theorem \re{mainTh5711} : partitions and graphs}\label{sec:Combinatorics}

The aim of this section is to introduce the combinatorial definitions which
will be used in the next section in order to  prove Theorem \re{mainTh5711}. The conventions we chose for the definitions of partitions, graphs, and for the less well-known notion of hypergraph are presented in 
Section \ref{appendicecombi} of the appendix.

\subsection{{Partitions}}\la{sect_partitions_19911}
\begin{itemize}
\item Let us recall that we denote by $\Part(k)$  the set of partitions of $V=\ensk$, and by  $\NC(k)$ the set of
  non-crossing partitions of $\ensk$ (a partition $\pi$ of $\ensk$  is said to be \emph{non-crossing} if there does not exist $x<y<z<t$ \st $x\simpi z\nsimpi y\simpi t$).
 \item For any given partition $\pi\in\Part(k)$, we denote by $\nc(\pi)$ be the minimal  non-crossing partition which is above
  $\pi$ for
  the refinement order;  the partitions induced by $\pi$ on the blocks of $\nc(\pi)$ are called the {\it connected components} of $\pi$. 
\item If $\nc(\pi)$ has only one block, then $\pi$ is said to be {\it connected}.
\ite  We define  $\op{thin}(\pi)$ to be the partition induced by $\pi$ on the subset of
$\ensk$ obtained by erasing $\ell$ whenever 
$
\ell+1\stackrel{\pi}{\sim} \ell
$ (with the convention $k+1=1$).
\ite A partition $\pi$ is said to be \emph{thin} if $\pi=\op{thin}(\pi)$.
\end{itemize}

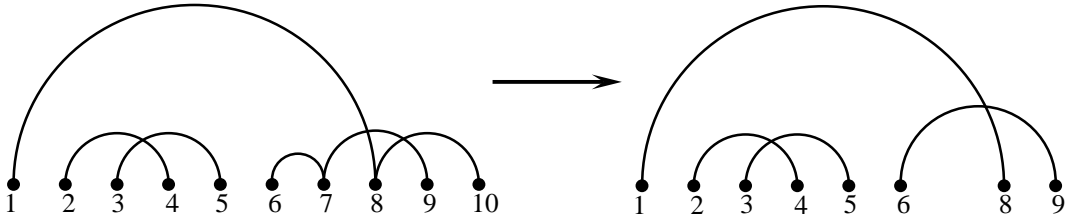
\begin{figure}[h!]
\scalebox{0.85} 
{
\begin{pspicture}(0,-0.4649999)(16.414999,2.8200002)
\psdots[dotsize=0.2](0.11411111,0.0)
\psdots[dotsize=0.2](0.91411114,0.0)
\psdots[dotsize=0.2](1.7141111,0.0)
\psdots[dotsize=0.2](2.514111,0.0)
\psdots[dotsize=0.2](3.314111,0.0)
\psdots[dotsize=0.2](4.114111,0.0)
\psdots[dotsize=0.2](4.914111,0.0)
\psdots[dotsize=0.2](5.7141113,0.0)
\psdots[dotsize=0.2](6.514111,0.0)
\psdots[dotsize=0.2](7.314111,0.0)
\usefont{T1}{ptm}{m}{n}
\rput(0.07473611,-0.2900001){1}
\usefont{T1}{ptm}{m}{n}
\rput(0.9698924,-0.2900001){2}
\usefont{T1}{ptm}{m}{n}
\rput(1.7370799,-0.2900001){3}
\usefont{T1}{ptm}{m}{n}
\rput(2.5769236,-0.2900001){4}
\usefont{T1}{ptm}{m}{n}
\rput(3.3427048,-0.2900001){5}
\usefont{T1}{ptm}{m}{n}
\rput(4.1600485,-0.2900001){6}
\usefont{T1}{ptm}{m}{n}
\rput(4.9600487,-0.2900001){7}
\usefont{T1}{ptm}{m}{n}
\rput(5.740361,-0.2900001){8}
\usefont{T1}{ptm}{m}{n}
\rput(6.55958,-0.2900001){9}
\usefont{T1}{ptm}{m}{n}
\rput(7.411611,-0.2900001){10}
\psarc[linewidth=0.04](2.914,0.0){2.8000002}{-1.5481567}{180.0}
\psarc[linewidth=0.04](6.514,0.0){0.8}{0.0}{180.0}
\psarc[linewidth=0.04](1.7141111,0.0){0.8}{0.0}{180.0}
\psarc[linewidth=0.04](2.514111,0.0){0.8}{0.0}{180.0}
\psarc[linewidth=0.04](5.7141113,0.0399999){0.8}{0.0}{180.0}
\psarc[linewidth=0.04](4.514111,0.0799999){0.4}{0.0}{180.0}
\psdots[dotsize=0.2](9.84,-0.0199999)
\psdots[dotsize=0.2](10.64,-0.0199999)
\psdots[dotsize=0.2](11.44,-0.0199999)
\psdots[dotsize=0.2](12.24,-0.0199999)
\psdots[dotsize=0.2](13.04,-0.0199999)
\psdots[dotsize=0.2](13.84,-0.0199999)
\psdots[dotsize=0.2](15.44,-0.0199999)
\psdots[dotsize=0.2](16.24,-0.0199999)
\usefont{T1}{ptm}{m}{n}
\rput(9.800625,-0.3099999){1}
\usefont{T1}{ptm}{m}{n}
\rput(10.695782,-0.3099999){2}
\usefont{T1}{ptm}{m}{n}
\rput(11.462969,-0.3099999){3}
\usefont{T1}{ptm}{m}{n}
\rput(12.302813,-0.3099999){4}
\usefont{T1}{ptm}{m}{n}
\rput(13.068594,-0.3099999){5}
\usefont{T1}{ptm}{m}{n}
\rput(13.885938,-0.3099999){6}
\usefont{T1}{ptm}{m}{n}
\rput(15.46625,-0.3099999){8}
\usefont{T1}{ptm}{m}{n}
\rput(16.28547,-0.3099999){9}
\psarc[linewidth=0.04](12.639889,-0.0199999){2.8000002}{-1.5481567}{180.0}
\psarc[linewidth=0.04](11.44,-0.0199999){0.8}{0.0}{180.0}
\psarc[linewidth=0.04](12.24,-0.0199999){0.8}{0.0}{180.0}
\psarc[linewidth=0.04](15.039889,0.0200001){1.2}{0.0}{180.0}
\psline[linewidth=0.06cm,arrowsize=0.05291667cm 3.0,arrowlength=2.0,arrowinset=0.4]{->}(7.525889,1.6000003)(9.525888,1.6000003)
\end{pspicture} 
}
\caption{Illustration of the operation $\pi\mapsto\op{thin}(\pi)$}\label{fig:thin-pi}
\end{figure}

For instance, for 
$\pi=\left\{\{1,8,10\},\{2,4\},\{3,5\},\{6,7,9\}\right\}\in\Part(10)$ the partition  $\nc(\pi)$ is   $\{\{1,6,7,8,9,10\},\{2,3,4,5\}\}$ and   Figure \ref{fig:thin-pi}  illustrates  of the operation $\pi\mapsto\op{thin}(\pi)$.

For any function $f$ defined on a set $E$, we denote by $\ker f $ the partition of $E$ whose blocks are the level sets of $f$. For any partition  $\pi\in \Part(k)$,   for each $\ell\in \ensk$, we denote by $\pi(\ell)$   the  index of the class
of $\ell$ in $\pi$,  after having ordered the classes according to the
order of their first element (for the example given Figure \ref{fig:pi}, we
have $\pi(1)=1$, $\pi({2})=2$, $\pi(3)=3$ and $\pi(4)=2$, $\pi(5)=3$, $\pi(6)=4)$,...).  Moreover, $0$ is identified with $k$ and $k+1$ is identified with $1$, so that $\pi(0)=\pi(k)$ and $\pi(k+1)=\pi(1)$.\\

\subsection{{Graphs}}

Let $G=(V,E=(e_i, i\in I))$ be a graph, and $\pi$ a partition of $V$. Let $\tilde{\pi}$ be the canonical surjection from $V$ onto $\pi$. The \emph{quotient graph} $G_{/\pi}$ is the graph with $\pi$ as set of vertices  and with edges $(e_i^{(\pi)},i\in I)$ such that $e_i^{(\pi)}$ is the edge between $\tilde{\pi}(v)$ and $\tilde{\pi}(w)$ if $e_i$ is an edge between $v$ and $w$, with same direction if the graph is directed. Note that the quotient graph can have strictly less vertices than the initial one, but it has as much edges as the initial one.
Note also that the quotient of  a circuit is still a circuit.

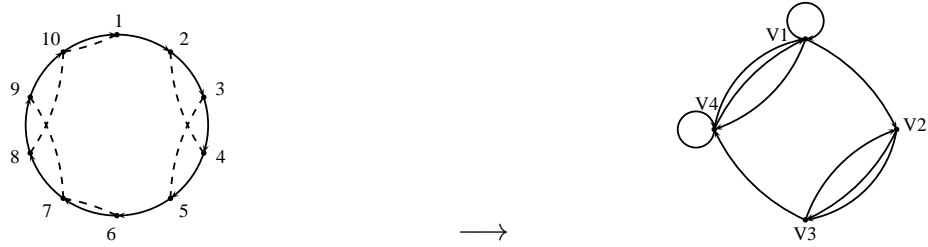
\begin{figure}[h!]
\centering
\subfigure[A cyclic graph $G$ and a partition $\pi$ of its vertices]{
\scalebox{0.6} 
{
\begin{pspicture}(0,-2.5443747)(11.5835,2.4834378)
\rput(-2.8,0){
\psdots[dotsize=0.12](9.176,2.0000002)
\psdots[dotsize=0.12](10.352,-1.6179997)
\psdots[dotsize=0.12](10.352,1.6180003)
\psdots[dotsize=0.12](11.078,0.6180001)
\psdots[dotsize=0.12](11.078,-0.6179997)
\psdots[dotsize=0.12](9.176,-1.9999998)
\psdots[dotsize=0.12](8.0,-1.6179997)
\psdots[dotsize=0.12](7.274,-0.6179997)
\psdots[dotsize=0.12](7.274,0.6180001)
\psdots[dotsize=0.12](8.0,1.6180003)
\usefont{T1}{ptm}{m}{n}
\rput(9.222875,2.3100002){1}
\usefont{T1}{ptm}{m}{n}
\rput(10.654593,1.9100002){2}
\usefont{T1}{ptm}{m}{n}
\rput(11.443656,0.81000024){3}
\usefont{T1}{ptm}{m}{n}
\rput(11.456938,-0.68999976){4}
\usefont{T1}{ptm}{m}{n}
\rput(10.645532,-1.8899997){5}
\usefont{T1}{ptm}{m}{n}
\rput(9.051312,-2.3899999){6}
\usefont{T1}{ptm}{m}{n}
\rput(7.6513124,-1.8899997){7}
\usefont{T1}{ptm}{m}{n}
\rput(6.94475,-0.68999976){8}
\usefont{T1}{ptm}{m}{n}
\rput(6.951156,0.81000024){9}
\usefont{T1}{ptm}{m}{n}
\rput(7.7285,1.9100002){10}
\psarc[linewidth=0.04,arrowsize=0.05291667cm 2.0,arrowlength=1.4,arrowinset=0.4]{<-}(9.176,0.0){2.0}{54.0}{90.0}
\psarc[linewidth=0.04,arrowsize=0.05291667cm 2.0,arrowlength=1.4,arrowinset=0.4]{<-}(9.176,0.0){2.0}{18.0}{54.0}
\psarc[linewidth=0.04,arrowsize=0.05291667cm 2.0,arrowlength=1.4,arrowinset=0.4]{<-}(9.176,0.0){2.0}{-18.0}{18.0}
\psarc[linewidth=0.04,arrowsize=0.05291667cm 2.0,arrowlength=1.4,arrowinset=0.4]{<-}(9.176,0.0){2.0}{-54.0}{-18.0}
\psarc[linewidth=0.04,arrowsize=0.05291667cm 2.0,arrowlength=1.4,arrowinset=0.4]{<-}(9.176,0.0){2.0}{-90.0}{-54.0}
\psarc[linewidth=0.04,arrowsize=0.05291667cm 2.0,arrowlength=1.4,arrowinset=0.4]{<-}(9.176,0.0){2.0}{90.0}{126.0}
\psarc[linewidth=0.04,arrowsize=0.05291667cm 2.0,arrowlength=1.4,arrowinset=0.4]{<-}(9.176,0.0){2.0}{126.0}{162.0}
\psarc[linewidth=0.04,arrowsize=0.05291667cm 2.0,arrowlength=1.4,arrowinset=0.4]{<-}(9.176,0.0){2.0}{162.0}{198.0}
\psarc[linewidth=0.04,arrowsize=0.05291667cm 2.0,arrowlength=1.4,arrowinset=0.4]{<-}(9.176,0.0){2.0}{198.0}{234.0}
\psarc[linewidth=0.04,arrowsize=0.05291667cm 2.0,arrowlength=1.4,arrowinset=0.4]{<-}(9.176,0.0){2.0}{234.0}{270.0}
\psarc[linewidth=0.04,linestyle=dashed,dash=0.16cm 0.16cm](14.352,1.6820003){4.0}{180.9}{215.1}
\psarc[linewidth=0.04,linestyle=dashed,dash=0.16cm 0.16cm](14.352,-1.6819999){4.0}{144.9}{179.1}
\psarc[linewidth=0.04,linestyle=dashed,dash=0.16cm 0.16cm](7.366,-5.5680003){4.0}{63.1}{80.9}
\psarc[linewidth=0.04,linestyle=dashed,dash=0.16cm 0.16cm](7.366,5.5680003){4.0}{-80.9}{-63.1}
\psarc[linewidth=0.04,linestyle=dashed,dash=0.16cm 0.16cm](4.0,1.6820004){4.0}{-35.1}{-0.9}
\psarc[linewidth=0.04,linestyle=dashed,dash=0.16cm 0.16cm](4.0,-1.6819996){4.0}{0.9}{35.1}
}
\end{pspicture} 
}\label{fig:pi-cycle}
}
\hspace{1cm}$\lto$\hspace{1cm}
\subfigure[The quotient $\Gspi$ of the left graph $G$ by the partition $\pi$]{
\scalebox{0.6} 
{
\begin{pspicture}(0,-2.4453125)(8.875,2.82)
\rput(-2,0){
\psdots[dotsize=0.12](6.23,2.0)
\psdots[dotsize=0.12](8.23,0.0)
\psdots[dotsize=0.12](6.23,-2.0)
\psdots[dotsize=0.12](4.23,0.0)
\usefont{T1}{ptm}{m}{n}
\rput(5.6164064,2.11){V1}
\usefont{T1}{ptm}{m}{n}
\rput(8.630938,0.11){V2}
\usefont{T1}{ptm}{m}{n}
\rput(6.2234373,-2.29){V3}
\usefont{T1}{ptm}{m}{n}
\rput(4.0721874,0.57){V4}
\psarc[linewidth=0.04,arrowsize=0.05291667cm 2.0,arrowlength=1.4,arrowinset=0.4]{<-}(4.584,-1.6459998){4.0}{24.3}{65.7}
\psarc[linewidth=0.04,arrowsize=0.05291667cm 2.0,arrowlength=1.4,arrowinset=0.4]{<-}(4.584,1.6459999){4.0}{-65.7}{-24.3}
\psarc[linewidth=0.04,arrowsize=0.05291667cm 2.0,arrowlength=1.4,arrowinset=0.4]{<-}(5.858,0.37200013){2.4}{-81.1}{-8.9}
\psarc[linewidth=0.04,arrowsize=0.05291667cm 2.0,arrowlength=1.4,arrowinset=0.4]{<-}(9.26,-3.03){3.2}{108.8}{161.2}
\psarc[linewidth=0.04,arrowsize=0.05291667cm 2.0,arrowlength=1.4,arrowinset=0.4]{<-}(7.876,1.6459999){4.0}{204.3}{245.7}
\psarc[linewidth=0.04,arrowsize=0.05291667cm 2.0,arrowlength=1.4,arrowinset=0.4]{<-}(7.876,-1.6459998){4.0}{114.3}{155.7}
\psarc[linewidth=0.04,arrowsize=0.05291667cm 2.0,arrowlength=1.4,arrowinset=0.4]{<-}(6.6019993,-0.37200013){2.4}{98.9}{171.1}
\psarc[linewidth=0.04,arrowsize=0.05291667cm 2.0,arrowlength=1.4,arrowinset=0.4]{<-}(3.2,3.03){3.2}{-71.2}{-18.8}
\psarc[linewidth=0.04,arrowsize=0.05291667cm 2.0,arrowlength=1.4,arrowinset=0.4]{<-}(6.23,2.4){0.4}{-90.0}{270.0}
\psarc[linewidth=0.04,arrowsize=0.05291667cm 2.0,arrowlength=1.4,arrowinset=0.4]{<-}(3.83,0.0){0.4}{0.0}{360.0}}
\end{pspicture} 
}\label{fig:vraigraphequotient}
}
\caption{A cyclic graph and one of its quotients}
\end{figure}

For instance, let us consider the preceding partition $\pi$ and $G$ be the cyclic graph with vertex set $V=\left\{1,\ldots,10\right\}$ and edges $1\to
2\to\cdots\to 10\to1$ (\emph{cf.} Figure \ref{fig:pi-cycle} where we  draw the partition $\pi$ with dashed lines).  The
quotient graph $G_{/\pi}$ is then given by Figure
\ref{fig:vraigraphequotient} with $V1=\{1,8,10\}$, $V2=\{2,4\}$, $V3=\{3,5\}$ and $V4=\{6,7,9\}$.

Let $G=(V,E)$ a directed graph. For any vertex $v\in
V$, we denote by $E^+(v)$ the subset of out-going edges from $v$, and by $E^-(v)$ the subset of in-going edges to
$v$:
\begin{eqnarray*}
E^+(v)&=&\left\{e\in E\ste(\exists v'\in \pi)(e=(v,v'))\right\}\\
E^-(v)&=&\left\{e\in E\ste (\exists v'\in \pi)(e=(v',v))\right\}
\end{eqnarray*}
Notice that $E^+(v)$ and $E^-(v)$  are not necessarily disjoints.\\

Let us consider $\{1,\ldots,N\}^E$ as a set of
``colorings'' of the edges of $G$ by the colors $1, \ld, N$. For any coloring $c\in\{1,\ldots,N\}^E$, let
$n_c^+(v,i):=\vert E^+(v)\cap c^{-1}(\{i\})\vert$ and $n_c^-(v,i):=\vert
E^-(v)\cap c^{-1}(\{i\})\vert$. If for each vertex $v$ and each color $i$, $n_c^+(v,i)=n_c^{-}(v,i)$, the coloring is called \emph{admissible}. \label{ncdefplus}  

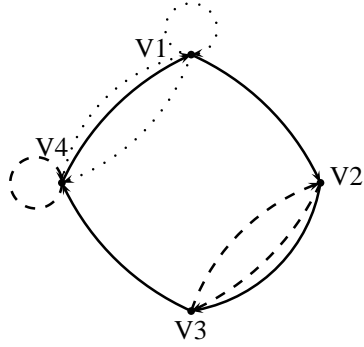
\begin{figure}[h!]
\scalebox{0.85} 
{
\begin{pspicture}(0,-2.4453125)(8.875,2.82)
\rput(-2,0){
\psdots[dotsize=0.12](6.23,2.0)
\psdots[dotsize=0.12](8.23,0.0)
\psdots[dotsize=0.12](6.23,-2.0)
\psdots[dotsize=0.12](4.23,0.0)
\usefont{T1}{ptm}{m}{n}
\rput(5.6164064,2.11){V1}
\usefont{T1}{ptm}{m}{n}
\rput(8.630938,0.11){V2}
\usefont{T1}{ptm}{m}{n}
\rput(6.2234373,-2.29){V3}
\usefont{T1}{ptm}{m}{n}
\rput(4.0721874,0.57){V4}
\psarc[linewidth=0.04,arrowsize=0.05291667cm 2.0,arrowlength=1.4,arrowinset=0.4]{<-}(4.584,-1.6459998){4.0}{24.3}{65.7}
\psarc[linewidth=0.04,linestyle=dashed,dash=0.16cm 0.16cm,arrowsize=0.05291667cm 2.0,arrowlength=1.4,arrowinset=0.4]{<-}(4.584,1.6459999){4.0}{-65.7}{-24.3}
\psarc[linewidth=0.04,arrowsize=0.05291667cm 2.0,arrowlength=1.4,arrowinset=0.4]{<-}(5.858,0.37200013){2.4}{-81.1}{-8.9}
\psarc[linewidth=0.04,linestyle=dashed,dash=0.16cm 0.16cm,arrowsize=0.05291667cm 2.0,arrowlength=1.4,arrowinset=0.4]{<-}(9.26,-3.03){3.2}{108.8}{161.2}
\psarc[linewidth=0.04,arrowsize=0.05291667cm 2.0,arrowlength=1.4,arrowinset=0.4]{<-}(7.876,1.6459999){4.0}{204.3}{245.7}
\psarc[linewidth=0.04,arrowsize=0.05291667cm 2.0,arrowlength=1.4,arrowinset=0.4]{<-}(7.876,-1.6459998){4.0}{114.3}{155.7}
\psarc[linewidth=0.04,linestyle=dotted,dotsep=0.16cm,arrowsize=0.05291667cm 2.0,arrowlength=1.4,arrowinset=0.4]{<-}(6.6019993,-0.37200013){2.4}{98.9}{171.1}
\psarc[linewidth=0.04,linestyle=dotted,dotsep=0.16cm,arrowsize=0.05291667cm 2.0,arrowlength=1.4,arrowinset=0.4]{<-}(3.2,3.03){3.2}{-71.2}{-18.8}
\psarc[linewidth=0.04,linestyle=dotted,dotsep=0.16cm,arrowsize=0.05291667cm 2.0,arrowlength=1.4,arrowinset=0.4]{<-}(6.23,2.4){0.4}{-90.0}{270.0}
\psarc[linewidth=0.04,linestyle=dashed,dash=0.16cm 0.16cm,arrowsize=0.05291667cm 2.0,arrowlength=1.4,arrowinset=0.4]{<-}(3.83,0.0){0.4}{0.0}{360.0}
}
\end{pspicture} 
}
\caption{An admissible coloring}\label{fig:graphequotient}
\end{figure}

Figure \ref{fig:graphequotient} presents  an instance of an admissible coloring of $G_{/\pi}$
with three colors.

A partition $\tau$ of the edges of $G$ will be said {\it admissible} if it is the kernel of an admissible coloring. In this case, for  any vertex $v\in G$, we define 
\be\la{571114h5}
n_\tau(v):=(n_\tau(v,1), n_\tau(v,2), \ld\ld)
\ee 
to be the common value  of the decreasing reordering   of the families $(n^+_c(v,1), n^+_c(v,2), \ld\ld)$ for $c$ coloring  \st $\ker c=\tau$.

\subsection{Hypergraph associated to a partition of the edges of a graph}\la{def_hypergraph_ass_color} Let $G$ be a graph with vertex set $V$, $\pi$ be a partition of $V$ and $\tau$ be a partition of the edge set of $G$. Then one can define $H(\pi, \tau)$ to be the hypergraph with the same vertex set as $G_{/\pi}$ (i.e. $\pi$) and with edges $(E_W, W\in \tau)$, where each edge $E_W$ is the set of blocks $J\in \pi$ \st at least one edge of $G_{/\pi}$ starting or ending at $J$ belongs to $W$.  

\begin{figure}[h!]
\scalebox{0.5} 
{
\begin{pspicture}(0,-3.279343)(6.46,3.1806571)
\rput(-3,0){
\psdots[dotsize=0.12](3.21,2.1906571)
\psdots[dotsize=0.12](5.21,0.19065711)
\psdots[dotsize=0.12](3.21,-1.8093429)
\psdots[dotsize=0.12](1.21,0.19065711)
\usefont{T1}{ptm}{m}{n}
\rput(2.5728126,2.300657){V1}
\usefont{T1}{ptm}{m}{n}
\rput(5.6018753,0.3006571){V2}
\usefont{T1}{ptm}{m}{n}
\rput(3.1868749,-2.0993428){V3}
\usefont{T1}{ptm}{m}{n}
\rput(1.044375,0.76065713){V4}
\rput{45.0}(1.5284188,-1.0486152){\psellipse[linewidth=0.04,linestyle=dotted,dotsep=0.16cm,dimen=outer](2.03,1.3206571)(2.0,0.77)}
\psellipse[linewidth=0.04,linestyle=dashed,dash=0.16cm 0.16cm,dimen=outer](3.18,-0.6893429)(2.6,1.8)
\pscircle[linewidth=0.04,dimen=outer](3.23,-0.049342893){3.23}
}
\end{pspicture} 
}
\caption{The graph $H(\pi,\tau)$ with $\pi$, $G_{/\pi}$ and $\tau$ as given by the preceding
example Figure \ref{fig:graphequotient}}\label{fig:hypergraph-color}
\end{figure}
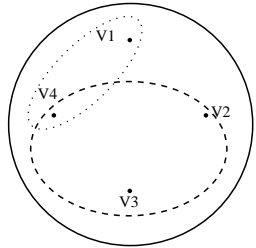

For instance, with $\pi$, $G_{/\pi}$ and $\tau$ as given by the preceding
example Figure \ref{fig:graphequotient},
  $H(\pi,\tau)$ is the hypergraph with three edges drawn in Figure \ref{fig:hypergraph-color}. Another example, where $H(\pi, \tau)$ has no cycle, is given at Figure \re{Correspondance_21911}.

 The next proposition will be used in the following. For $\tau$ a partition of the edges of a graph, a cycle of the graph is said to be $\tau$-\emph{monochromatic} if all the edges it visits belong to the same block of $\tau$.

\beg{propo}\la{prop_fond_real_14911}Let $G=(V,E)$ be a   directed graph, which is a circuit.   Fix  $\pi\in \Part(V)$ and  $\tau$  be a partition of the edge set of $G$ \st $H(\pi,\tau)$ has no cycle.  Then $\Gspi$ is a disjoint union of $\tau$-monochromatic cycles. As a consequence, $\tau$ is admissible.\en{propo}

\beg{pr}If  $\tau$ has more than one block (the case with one bloc being obvious), since $\Gspi$ is a circuit, one can find a closed path 
 $$\gamma=(v_0,e_{1},v_1,e_{2},\ldots,e_{k},v_{k})\qquad\trm{  (with $v_k=v_0$)}$$ of $\Gspi$ which visits each edge of $\Gspi$ exactly once and \st $e_1$ and $e_k$ do not belong to the same block of $\tau$.  
For each edge $e$ of $\Gspi$,  let $E(e)$ be the edge of $H(\pi, \tau)$  consisting of all vertices of $\Gspi$ which are the beginning or the end of an edge of $\Gspi$ with the same  color as $e$ (i.e. in the same block of $\tau$ as $e$).  Then $(v_0,E(e_{1}),v_1,E(e_{2}),\ldots,E(e_{k}),v_{k})$ is a path of $H(\pi, \tau)$, i.e. for all $\ell=1, \le, k$, $v_{\ell-1}, v_\ell\in E(e_\ell)$. 
Let us introduce $i_0=0<i_1<\cd <i_p=k$ \st  $i_1, \ld ,i_p$ are the instants  where the color of the edges change in $\gamma$ (i.e. for all $s=1, \ld, p$, $e_{i_{s-1}+1}, \ld, e_{i_s}$ have the same color, and for all $s=1, \ld,p-1$, $e_{i_s}$  and $e_{i_s+1}$ do not have the same color. Then $$\bar{\gamma}
:=(v_{i_0}, E(e_{i_1}), v_{i_1}, \ld\ld, v_{i_{p-1}}, E(e_{i_p}), v_{i_p})$$ is a also a path of $H(\pi, \tau)$.

Since $H(\pi,\tau)$ has no cycle, hence is linear, $v_{i_r}=E(e_{i_r})\cap E(e_{i_{r+1}})$ for $r=1,\ldots,p-1$ and $v_{i_0}=v_{i_p}=E(e_{i_1})\cap E(e_{i_p})$. Let $v_{i_r},v_{i_{r+1}},\ldots,v_{i_s}$ a sequence of pairwise distinct vertices, with $v_{i_r}=v_{i_s}$ and $r<s$ (this exists since $v_{i_0}=v_{i_p}$). Then $s=r+1$ otherwise $\left( v_{i_r},E(e_{i_{r+1}}),v_{i_{r+1}},\ldots,E(e_{i_s}),v_{i_s}\right)$ would be a cycle of $H(\pi,\tau)$. Let us consider the path $(v_{i_r},e_{i_r+1},v_{i_r+1},\ldots,e_{i_s},v_{i_s})$. All the edges have the same color, and this is a disjoint union of $\tau$-monochromatic cycles. 

Let us now  apply the same trick to the reduced path
$$
(v_{i_0}, E(e_{i_1}), v_{i_1}, \ld,v_{i_r},E(e_{i_{s+1}}),\ld, v_{i_{p-1}}, E(e_{i_p}), v_{i_p}).
$$

This allows us to achieve by induction the proof of the proposition.

\en{pr}

An example of $\tau$ \st $H(\pi, \tau)$ has no cycle is given at  Figure \re{Correspondance_21911} below (another example of hypergraph with no cycle is given in Figure \re{fig:hypergraph_with_no_cycle} of the appendix).

\begin{center}
\begin{figure}[h!]
\scalebox{0.5} 
{
\begin{pspicture}(0,-2.72)(19.02,3.28)
\psdots[dotsize=0.12](16.33,2.31)
\psdots[dotsize=0.12](18.33,0.31)
\psdots[dotsize=0.12](16.33,-1.69)
\psdots[dotsize=0.12](14.33,0.31)
\usefont{T1}{ptm}{m}{n}
\rput(15.669219,2.4199998){V1}
\usefont{T1}{ptm}{m}{n}
\rput(18.712812,0.42){V2}
\usefont{T1}{ptm}{m}{n}
\rput(16.290312,-1.9799999){V3}
\usefont{T1}{ptm}{m}{n}
\rput(14.156563,0.88){V4}
\rput{45.0}(5.455566,-10.290901){\psellipse[linewidth=0.04,linestyle=dotted,dotsep=0.16cm,dimen=outer](15.15,1.44)(2.0,0.77)}
\psdots[dotsize=0.12](6.23,2.0)
\psdots[dotsize=0.12](8.23,0.0)
\psdots[dotsize=0.12](6.23,-2.0)
\psdots[dotsize=0.12](4.23,0.0)
\usefont{T1}{ptm}{m}{n}
\rput(5.5928125,2.11){V1}
\usefont{T1}{ptm}{m}{n}
\rput(8.621876,0.11){V2}
\usefont{T1}{ptm}{m}{n}
\rput(6.206875,-2.29){V3}
\usefont{T1}{ptm}{m}{n}
\rput(4.064375,0.57){V4}
\psarc[linewidth=0.04,arrowsize=0.05291667cm 2.0,arrowlength=1.4,arrowinset=0.4]{<-}(4.584,-1.6459998){4.0}{24.3}{65.7}
\psarc[linewidth=0.04,arrowsize=0.05291667cm 2.0,arrowlength=1.4,arrowinset=0.4]{<-}(4.584,1.6459998){4.0}{-65.7}{-24.3}
\psarc[linewidth=0.04,arrowsize=0.05291667cm 2.0,arrowlength=1.4,arrowinset=0.4]{<-}(5.858,0.37200013){2.4}{-81.1}{-8.9}
\psarc[linewidth=0.04,arrowsize=0.05291667cm 2.0,arrowlength=1.4,arrowinset=0.4]{<-}(9.26,-3.03){3.2}{108.8}{161.2}
\psarc[linewidth=0.04,arrowsize=0.05291667cm 2.0,arrowlength=1.4,arrowinset=0.4]{<-}(9.66,0.0){4.0}{150.0}{210.0}
\psarc[linewidth=0.04,linestyle=dotted,dotsep=0.16cm,arrowsize=0.05291667cm 2.0,arrowlength=1.4,arrowinset=0.4]{<-}(6.6019993,-0.37200013){2.4}{98.9}{171.1}
\psarc[linewidth=0.04,linestyle=dotted,dotsep=0.16cm,arrowsize=0.05291667cm 2.0,arrowlength=1.4,arrowinset=0.4]{<-}(3.2,3.03){3.2}{-71.2}{-18.8}
\psline[linewidth=0.06cm,arrowsize=0.05291667cm 3.0,arrowlength=2.0,arrowinset=0.4]{->}(10.02,0.08)(12.42,0.08)
\rput{-90.0}(16.94,17.5){\psellipse[linewidth=0.04,dimen=outer](17.22,0.28)(3.0,1.8)}
\end{pspicture} 
}\caption{On the left: an admissible partition $\tau$ of the edges of $G_{/\pi}$ whose corresponding hypergraph $H(\pi, \tau)$ (on the right) has no cycle.}\la{Correspondance_21911}
\end{figure}
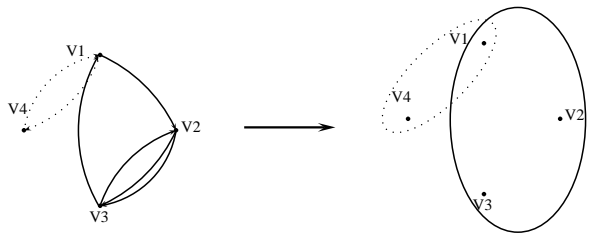
\en{center}

\section{Proof of Theorem \re{mainTh5711}}\la{771113h37}

The proof of Theorem \re{mainTh5711} will follow from Propositions
\re{571119h07} and \re{571119h24} and Lemma \re{lem:0519-1} below.   We suppose that Hypothesis  \re{hyp_NfixeUn}  holds. 

\subsection{Convergence   of the mean spectral distribution in the case with moments}
In this section, we prove the following proposition. The parametrization of infinitely divisible laws via  formula $\mu=\nu_*^{\gamma,\si}$ is introduced in Section \re{appendiceinfdiv} of the appendix.

\begin{propo}\la{571119h07}
Consider an infinitely divisible law
  $\mu=\nu_*^{\gamma,\si}$ \st $\si$ has compact support. Define the mean
  spectral distribution $\bar{\Lam}_N^{(\mu)}$ of a $\PUNmu$-distributed
  random matrix $M$ by 
$$
\qquad\qquad\int g(t)\bar{\Lam}_N^{(\mu)}(\ud t)=\E\left[\ff{N}\Tr
g(M)\right]\qquad\qquad\trm{ (for any Borel function $g$).}
$$   
Then $\bar{\Lam}_N^{(\mu)}$ converges weakly, as $N\lto\infty$,  to a \pro
measure $\Lam_\Gamma(\mu)$ which depends only on $\mu$ and on the function
$\Gamma$ introduced in Hypothesis  \re{hyp_NfixeUn}. This distribution admits exponential moments of any order, hence is
characterized by its moments, given by Formula \eqre{eq:0519-2ter} in
Theorem \ref{mainTh5711}.
\end{propo}

 Let us introduce a $\PUNmu$-distributed random matrix $M$ and fix $k\ge 1$.  Proposition \re{571119h07} follows directly from Lemmas \re{571119h10},    \re{pte:0213-3gusgus} and \re{lem_carleman} below.

\begin{lem}\la{571119h10}
For any $N\ge 1$ and $k\geq 1$,  there is a collection $(f_N(\pi))_{\pi\in \Part(k)}$ of
numbers depending only on the distribution of $U_N$ \st 
$$
\E\left[ \ff{N}\Tr M^k\right]=\sum_{\pi\in \Part(k)}f_N(\pi) c_\pi(\mu).
$$
Moreover, for any $\pi$, $f_N(\pi)$ is given by the formula 
\be\la{f_Nformule1}
f_N(\pi)=N^{\Card{\pi} -1-k }\,\E\left[\prod_{\ell=1}^k
U_N^{\pi(\ell)*}U_N^{\pi(\ell+1)}\right],
\ee 
where    $(U_N^{i})_{i\ge 1}$ is a sequence of i.i.d.  copies of $U_N$.  
\end{lem} 

\beg{pr}
For each $n\ge 1$, set $\nu_n=\mu^{*\ff{n}}$ and let $X_n^1, X_n^2,
\ld $ be a sequence of i.i.d. $\nu_n$-distributed random variables,
independent of the $U_N^i$'s, so that according to Theorem \re{defPUNmu}, for $$
M_n:=\ff{N}\sum_{i=1}^{N\ti n}X_n^{i}\cdot U_N^{i}U_N^{i*},
$$ 
the law of $M_n$ converges weakly to $
\PUNmu.
$ 
Then it is immediate that for $n$ large enough, 
$$
\E\left[\ff{N}\Tr M_n^k\right]=\ff{N^{k+1}}
\sum_{\pi\in \Part(k)}\frac{(Nn)!}{(Nn-\Card{\pi})!}  \,m_\pi(\nu_n)\,\E\left[\Tr \prod_{\ell=1}^k
  U_N^{\pi(\ell)}U_N^{\pi(\ell)*}\right],
$$ 
where for any $\pi$, 
$$
m_\pi(\nu_n):=\E\left[\prod_{\ell=1}^k X_n^{\pi(\ell)}\right]=\prod_{J\in \pi}
\int t^{\Card{J}}\ud \nu_n(t).
$$ 
Note that since the $U_N^i$'s are column vectors, $\Tr \prod_{\ell=1}^k U_N^{\pi(\ell)}U_N^{\pi(\ell)*}=\prod_{\ell=1}^k U_N^{\pi(\ell)*}U_N^{\pi(\ell+1)}$. Moreover,  by \cite[Th. 1.6]{FBDAOPinfdiv}, 
$
n^{\Card{\pi}}m_\pi(\nu_n)\ninf c_\pi(\mu),
$ 
which allows to conclude.
\end{pr}



Let us now consider the large $N$ limit of $f_N(\pi)$. Let $G$ be the cyclic graph with vertex set $V=\{1,\ldots,k\}$ and edges
$1\rightarrow2\rightarrow\cdots\rightarrow k\rightarrow 1$.  The hypergraph $H(\pi,\tau)$ is defined in Section \re{def_hypergraph_ass_color} and the operator $\op{thin}$ on partitions is  defined in Section \re{sect_partitions_19911}.

\beg{lem}\la{pte:0213-3gusgus}
For each $\pi\in \Part(k)$, there is $f_\Ga(\pi)\ge 0$ \st  
$
f_N(\pi)\Ninf  f_\Ga(\pi)$. The number $f(\pi)$ factorizes along the connected components of $\pi$ and is given by:\be\la{20911.19h}f_\Ga(\pi)=\sum_{\substack{\tau\trm{ partition of the edges of}\\  \trm{$\Gspi$ s.t.  $H(\pi,\tau)$ has no cycle}}}\prod_{J\in
  \pi}\Gamma(n_\tau(J)).\ee
  If, moreover, $U_N$ is a unitary vector,   for any $\pi$, $f_\Ga(\pi)=f_\Ga(\op{thin}(\pi))$ and for
  $\pi$   non crossing    $f_\Ga(\pi)=1$.
\en{lem}

Before proving the lemma, let us give two examples to make the set of partitions $\tau$ \st $H(\pi, \tau)$ has no cycle   a bit more intuitive. 

a) Consider for example the case $k=4$ and $\pi=\{\{1\},\{2\},\{3\},\{4\}\}$ (trivial partition). Then it can easily be seen that there is only one such $\tau$: it is the trivial partition with only one block.  

b)  Consider now the partition $\pi=\{\{1\},\{2,4\},\{3\}\}$. Then there are two such $\tau$'s: the first one is the one with only one block, containing all vertices;   the second one is  $\{\{(1,2),(4,1)\},\{(2,3),(3,4)\}\}$. 
  
\beg{pr}$\bullet$  Let us first prove that $f_N(\pi)$ has a limit $f_\Ga(\pi)$ given by Formula \eqre{20911.19h} as $N\to\infty$. We start with Formula \eqre{f_Nformule1}, expand the matrix product, and use the independence of the $U_N^i$'s and the exchangeability of the entries $U_N(1),\ld, U_N(N)$ of $U_N$: 
\beq
f_N(\pi)&=&N^{\Card{\pi} -1 -k}\,\E\left[\prod_{\ell=1}^k
U_N^{\pi(\ell)*}U_N^{\pi(\ell+1)}\right]\\
&=&N^{\Card{\pi} -1-k }\sum_{1\le i_1, \ld, i_k\le N}\E\left[\prod_{\ell=1}^k
\ovl{U_N^{\pi(\ell)}(i_\ell)}U_N^{\pi(\ell+1)}(i_\ell)\right]\\
&=&N^{\Card{\pi} -1 -k}\sum_{1\le i_1, \ld, i_k\le N}\prod_{J\in \pi}\E\left[\prod_{\ell\in J}
\ovl{U_N(i_\ell)}\prod_{\ell'\trm{ s.t. } \ell'+1\in J}U_N(i_{\ell'})\right]\\
&=&N^{\Card{\pi} -1-k }\sum_{\substack{\tau\trm{ partition}\\ \trm{of the edges of $\Gspi$}}}A_N^{|\tau|}\prod_{J\in \pi}\E\left[\prod_{\ell\in J}
\ovl{U_N(i_l^\tau)}\prod_{\ell'\trm{ s.t. } \ell'+1\in J}U_N(i_{l'}^{\tau})\right],
\eeq
where for all $\tau$, $A_N^{|\tau|}=N(N-1)\cd (N-|\tau|+1)$ and $i^\tau$ is a coloring of the edges of $\Gspi$ with kernel $\tau$. 

Now, for each $\tau$, for each $J$, let $a_\tau(J)$ be the number of blocks of $\tau$ containing an edge having $J$ as an extremity. By Formula \eqre{571117h30_mjj} of Hypothesis \ref{hyp_NfixeUn}, as $N\to\infty$,  the term associated to $\tau$ is $$O(N^{\vert\pi\vert-1+\vert\tau\vert-\sum_{J\in\pi}a_\tau(J)}),$$
and that in the particular case where $\tau$ is admissible, this term is  equivalent to
$$
N^{\vert\pi\vert+\vert\tau\vert-1-\sum_{J\in\pi}a_\tau(J)}\prod_{J\in
  \pi}\Gamma(n_\tau(J)).
$$

Let us consider the hypergraph $H(\pi,\tau)$ defined at Section \re{def_hypergraph_ass_color}, with edges $\left(E_W,W\in\tau\right)$. 
Notice   that each edge of
$G_{/\pi}$ is included in one edge of $H(\pi,\tau)$. This implies that $H(\pi,\tau)$ is  connected. Notice also that:
$$
\sum_{J\in \pi}a_\tau(J)=\sum_{J\in \pi \atop W\in\tau}\one_{J\in E_W}=\sum_{W\in\tau}\vert E_W\vert
$$
Hence, using Property \ref{prop:0503}, we know that 
$$
\vert\pi\vert+\vert\tau\vert-1-\sum_{J\in\pi}a(J)=\vert\pi\vert+\vert\tau\vert-1-\sum_{W\in\tau}\vert E_W\vert
$$
is non positive, and vanishes if and only if $H(\pi,\tau)$ has no cycle. By Proposition \re{prop_fond_real_14911}, this proves that 
$f_N(\pi)$ converges when $N$ goes to infinity, with a limit given by \eqre{20911.19h}.

$\bullet$ Let us now prove that $f_\Ga(\pi)$ factorizes along the connected components of $\pi$. 

Suppose that $\pi$ has at least two connected components, and let   $G$ be the cyclic graph on $\{1,\ldots,k\}$ introduced in the statement of the previous lemma  and $G_{/\pi}=(\pi,E)$ be the associated quotient graph. There exists at least one connected component $\pi_{\vert A}$ of $\pi$  whose support
$A$ is an interval $\left[l,l'\right]$. The corresponding subgraph of $G_{/\pi}$ is
connected to the rest of the graph through two directed edges, $\varepsilon=(V_1,V_2)$ and $\varepsilon'=(V_3,V_4)$ with $l-1\in V_1$, $l\in V_2$, $l'\in V_3$ and $l'+1\in V_4$ (with the usual identifications   $0=k$ and $1=k+1$).

For instance, if we still consider the partition
$\pi=\left\{\{1,8,10\},\{2,4\},\{3,5\},\{6,7,9\}\right\}\in\Part(10)$, then the restriction of $\pi$ to 
the interval $\{2,3,4,5\}$ is a connected component of $\pi$, and 
the corresponding subgraph of $G_{/\pi}$ (\emph{cf.} Figure \ref{fig:graphequotient} page \pageref{fig:graphequotient}) is its restriction to the set of vertices
$\{V2,V3\}$, and is connected to the rest of the graph by the edges $\varepsilon=(V1,V2)$ and
$\varepsilon'=(V3,V4)$:

\begin{figure}[h!]
\scalebox{0.8} 
{
\begin{pspicture}(0,-2.4453125)(8.875,2.82)
\rput(-1,0){
\psdots[dotsize=0.12](6.23,2.0)
\psdots[dotsize=0.12](8.23,0.0)
\psdots[dotsize=0.12](6.23,-2.0)
\psdots[dotsize=0.12](4.23,0.0)
\usefont{T1}{ptm}{m}{n}
\rput(5.6164064,2.11){V1}
\usefont{T1}{ptm}{m}{n}
\rput(8.630938,0.11){V2}
\usefont{T1}{ptm}{m}{n}
\rput(6.2234373,-2.29){V3}
\usefont{T1}{ptm}{m}{n}
\rput(4.0721874,0.57){V4}
\psarc[linewidth=0.04,arrowsize=0.05291667cm 2.0,arrowlength=1.4,arrowinset=0.4]{<-}(4.584,-1.6459998){4.0}{24.3}{65.7}
\psarc[linewidth=0.04,arrowsize=0.05291667cm 2.0,arrowlength=1.4,arrowinset=0.4]{<-}(4.584,1.6459999){4.0}{-65.7}{-24.3}
\psarc[linewidth=0.04,arrowsize=0.05291667cm 2.0,arrowlength=1.4,arrowinset=0.4]{<-}(5.858,0.37200013){2.4}{-81.1}{-8.9}
\psarc[linewidth=0.04,arrowsize=0.05291667cm 2.0,arrowlength=1.4,arrowinset=0.4]{<-}(9.26,-3.03){3.2}{108.8}{161.2}
\psarc[linewidth=0.04,arrowsize=0.05291667cm 2.0,arrowlength=1.4,arrowinset=0.4]{<-}(7.876,1.6459999){4.0}{204.3}{245.7}
\psarc[linewidth=0.04,arrowsize=0.05291667cm 2.0,arrowlength=1.4,arrowinset=0.4]{<-}(7.876,-1.6459998){4.0}{114.3}{155.7}
\psarc[linewidth=0.04,arrowsize=0.05291667cm 2.0,arrowlength=1.4,arrowinset=0.4]{<-}(6.6019993,-0.37200013){2.4}{98.9}{171.1}
\psarc[linewidth=0.04,arrowsize=0.05291667cm 2.0,arrowlength=1.4,arrowinset=0.4]{<-}(3.2,3.03){3.2}{-71.2}{-18.8}
\psarc[linewidth=0.04,arrowsize=0.05291667cm 2.0,arrowlength=1.4,arrowinset=0.4]{<-}(6.23,2.4){0.4}{-90.0}{270.0}
\psarc[linewidth=0.04,arrowsize=0.05291667cm 2.0,arrowlength=1.4,arrowinset=0.4]{<-}(3.83,0.0){0.4}{0.0}{360.0}
\usefont{T1}{ptm}{m}{n}
\rput(7.548125,1.41){$\varepsilon$}
\usefont{T1}{ptm}{m}{n}
\rput(4.918125,-1.41){$\varepsilon'$}
}
\end{pspicture} 
}
\caption{The quotient graph $G_{/\pi}$, with the edges $\eps$ and $\eps'$}\label{fig:graphequotientepsilon}
\end{figure}
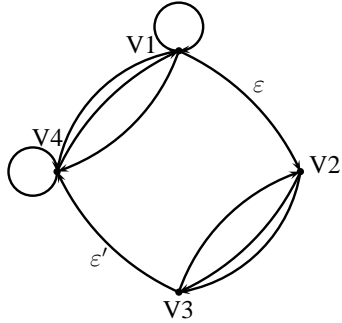

 Let us define the graph $T\left(G_{/\pi}\right)$ obtained from $G_{/\pi}$ by replacing the edges $\eps$ and $\eps'$ by the edges $\varepsilon_0=(V_3,V_2)$ and $\varepsilon_1=(V_1,V_4)$. This graph has two connected components, one on the subset of vertices $\pi_{\vert A}$, the other on the complementary $\pi_{\vert A^c}$. The first (resp. second) one, that we denote by $G_1$ (resp. $G_0$), is the restriction of $G_{/\pi}$ to $A$ (resp. $A^c$), plus the edge $\eps_1$ (resp. $\eps_0$).
 
For instance, $T(G_{/\pi})$ is drawn for the preceding example in Figure \ref{fig:t-graphequotient}.

 Let us consider a partition $\tau$ of the edges of $\Gspi$ \st $H(\pi, \tau)$ has no cycle. Since, by Proposition \re{prop_fond_real_14911}, $\Gspi$ is a union of pairwise disjoint $\tau$-monochromatic cycles,  $\eps$ and $\eps'$ belong to the same block of $\tau$.   Let us   define the partition $T(\tau)$ of $T\left(G_{/\pi}\right)$ deduced from $\tau$ by replacing $\eps$ and $\eps'$ by $\eps_0$ and $\eps_1$ in their block of $\tau$ (see Figure \re{fig:coloring-t-G}). It is easy to see that $T$ is a bijection between the set of partitions $\tau$ of the edges of $\Gspi$ \st $H(\pi, \tau)$ has no cycle and the set of pairs of such partitions of the two connected components  of $T(\Gspi)$. Moreover, $n_\tau(J)=n_{T(\tau)}(J)$ for any $J\in\pi$. The conclusion follows directly.

\begin{figure}[h!]\centering
\subfigure[The graph $T(G_{/\pi})$]{
\scalebox{0.8} 
{
\begin{pspicture}(0,-2.4453125)(8.865938,2.82)
\psdots[dotsize=0.12](6.23,2.0)
\psdots[dotsize=0.12](8.23,0.0)
\psdots[dotsize=0.12](6.23,-2.0)
\psdots[dotsize=0.12](4.23,0.0)
\usefont{T1}{ptm}{m}{n}
\rput(5.5928125,2.11){V1}
\usefont{T1}{ptm}{m}{n}
\rput(8.621876,0.11){V2}
\usefont{T1}{ptm}{m}{n}
\rput(6.206875,-2.29){V3}
\usefont{T1}{ptm}{m}{n}
\rput(4.064375,0.57){V4}
\psarc[linewidth=0.04,arrowsize=0.05291667cm 2.0,arrowlength=1.4,arrowinset=0.4]{<-}(4.584,1.6459999){4.0}{-65.7}{-24.3}
\psarc[linewidth=0.04,arrowsize=0.05291667cm 2.0,arrowlength=1.4,arrowinset=0.4]{<-}(5.858,0.37200013){2.4}{-81.1}{-8.9}
\psarc[linewidth=0.04,arrowsize=0.05291667cm 2.0,arrowlength=1.4,arrowinset=0.4]{<-}(9.26,-3.03){3.2}{108.8}{161.2}
\psarc[linewidth=0.04,arrowsize=0.05291667cm 2.0,arrowlength=1.4,arrowinset=0.4]{<-}(7.876,-1.6459998){4.0}{114.3}{155.7}
\psarc[linewidth=0.04,arrowsize=0.05291667cm 2.0,arrowlength=1.4,arrowinset=0.4]{<-}(6.6019993,-0.37200013){2.4}{98.9}{171.1}
\psarc[linewidth=0.04,arrowsize=0.05291667cm 2.0,arrowlength=1.4,arrowinset=0.4]{<-}(3.2,3.03){3.2}{-71.2}{-18.8}
\psarc[linewidth=0.04,arrowsize=0.05291667cm 2.0,arrowlength=1.4,arrowinset=0.4]{<-}(6.23,2.4){0.4}{-90.0}{270.0}
\psarc[linewidth=0.04,arrowsize=0.05291667cm 2.0,arrowlength=1.4,arrowinset=0.4]{<-}(3.83,0.0){0.4}{0.0}{360.0}
\psarc[linewidth=0.04,arrowsize=0.05291667cm 2.0,arrowlength=1.4,arrowinset=0.4]{<-}(8.23,-2.0){2.0}{90.0}{180.0}
\psarc[linewidth=0.04,arrowsize=0.05291667cm 2.0,arrowlength=1.4,arrowinset=0.4]{<-}(4.23,2.0){2.0}{-90.0}{0.0}
\usefont{T1}{ptm}{m}{n}
\rput(6.956875,-0.13){$\varepsilon_0$}
\usefont{T1}{ptm}{m}{n}
\rput(6.08125,0.73){$\varepsilon_1$}
\end{pspicture} 
}\label{fig:t-graphequotient}
}
\subfigure[The coloring $T(\tau)$]{
\scalebox{0.8} 
{
\begin{pspicture}(0,-2.4453125)(8.865938,2.82)
\psdots[dotsize=0.12](6.23,2.0)
\psdots[dotsize=0.12](8.23,0.0)
\psdots[dotsize=0.12](6.23,-2.0)
\psdots[dotsize=0.12](4.23,0.0)
\usefont{T1}{ptm}{m}{n}
\rput(5.5928125,2.11){V1}
\usefont{T1}{ptm}{m}{n}
\rput(8.621876,0.11){V2}
\usefont{T1}{ptm}{m}{n}
\rput(6.206875,-2.29){V3}
\usefont{T1}{ptm}{m}{n}
\rput(4.064375,0.57){V4}
\psarc[linewidth=0.04,linestyle=dashed,dash=0.16cm 0.16cm,arrowsize=0.05291667cm 2.0,arrowlength=1.4,arrowinset=0.4]{<-}(4.584,1.6459999){4.0}{-65.7}{-24.3}
\psarc[linewidth=0.04,arrowsize=0.05291667cm 2.0,arrowlength=1.4,arrowinset=0.4]{<-}(5.858,0.37200013){2.4}{-81.1}{-8.9}
\psarc[linewidth=0.04,linestyle=dashed,dash=0.16cm 0.16cm,arrowsize=0.05291667cm 2.0,arrowlength=1.4,arrowinset=0.4]{<-}(9.26,-3.03){3.2}{108.8}{161.2}
\psarc[linewidth=0.04,arrowsize=0.05291667cm 2.0,arrowlength=1.4,arrowinset=0.4]{<-}(7.876,-1.6459998){4.0}{114.3}{155.7}
\psarc[linewidth=0.04,linestyle=dotted,dotsep=0.16cm,arrowsize=0.05291667cm 2.0,arrowlength=1.4,arrowinset=0.4]{<-}(6.6019993,-0.37200013){2.4}{98.9}{171.1}
\psarc[linewidth=0.04,linestyle=dotted,dotsep=0.16cm,arrowsize=0.05291667cm 2.0,arrowlength=1.4,arrowinset=0.4]{<-}(3.2,3.03){3.2}{-71.2}{-18.8}
\psarc[linewidth=0.04,linestyle=dotted,dotsep=0.16cm,arrowsize=0.05291667cm 2.0,arrowlength=1.4,arrowinset=0.4]{<-}(6.23,2.4){0.4}{-90.0}{270.0}
\psarc[linewidth=0.04,linestyle=dashed,dash=0.16cm 0.16cm,arrowsize=0.05291667cm 2.0,arrowlength=1.4,arrowinset=0.4]{<-}(3.83,0.0){0.4}{0.0}{360.0}
\psarc[linewidth=0.04,arrowsize=0.05291667cm 2.0,arrowlength=1.4,arrowinset=0.4]{<-}(8.23,-2.0){2.0}{90.0}{180.0}
\psarc[linewidth=0.04,arrowsize=0.05291667cm 2.0,arrowlength=1.4,arrowinset=0.4]{<-}(4.23,2.0){2.0}{-90.0}{0.0}
\usefont{T1}{ptm}{m}{n}
\rput(6.956875,-0.13){$\varepsilon_0$}
\usefont{T1}{ptm}{m}{n}
\rput(6.08125,0.73){$\varepsilon_1$}
\end{pspicture} 
}\label{fig:coloring-t-G}
}
\caption{The graph $T(G_{/\pi})$ and the coloring $T(\tau)$} 
\end{figure}
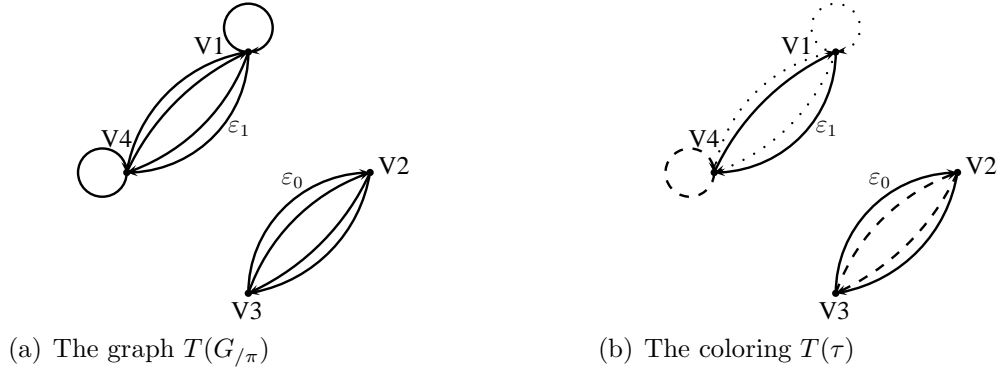

$\bullet$ Under the additional  hypothesis that $U_N$ is a unitary vector, it  follows directly from Formula \eqre{f_Nformule1} that $f_N(\pi)=f_N(\op{thin}(\pi))$, hence $f_\Ga(\pi)=f_\Ga(\op{thin}(\pi))$. It also follows from Formula \eqre{f_Nformule1} that for $k=1$, $f_N(\{\{1\}\})=1$. If $\pi$ is non crossing, then the partition induced by $\pi$ on its connected components are all some one-bloc partitions. Hence the fact that $f_\Ga$ factorizes on its connected components and  the identities $f_\Ga(\pi)=f_\Ga(\op{thin}(\pi))$ and $f_\Ga(\{\{1\}\})=1$ iterated imply that $f_\Ga(\pi)=1$.
\en{pr}

\beg{lem}\la{lem_carleman}Let us define $m_0=1$ and, for each $k\ge 1$, $$m_k=\sum_{\pi\in \Part(k)}f_\Ga(\pi) c_\pi(\mu).$$ Then the radius of convergence of the series $\sum_{k\ge 0} \f{m_k}{k!}z^k$ is infinite. 
\end{lem}

\beg{pr}{\bf Claim:} it suffices to prove the result under the additional hypothesis that $\Ga(\cdot)=1$.

Let us prove the claim.
Recall that $\gamma, \si$ are \st $\mu=\nu_*^{\gamma,\si}$. Let $\tilde{\si}$ be the push-forward of the measure $\si$ by the map $t\mapsto|t|$ and define 
 $\tmu:=\nu_*^{|\gamma|,\tilde{\si}}$.  Note that by Formula \eqre{571119h57} of the appendix, for any $k$, $|c_k(\mu)|\le c_k(\tmu)$.   
 Let us define, the function $\tGa(\cdot)$ by $\tGa(n_1, \ld, n_k):=1$ for any integers $n_1, \ld, n_k$.  For $C$ the constant of  Equation \eqre{571117h30ccbain}, we have, for any $k\ge 1$ and any $\pi\in \Part(k)$, $$0\le f_\Ga(\pi)\le C^k f_{\tGa}(\pi).$$
 As a consequence, if one defines, for each $k\ge 0$, $$\tilde{m}_k:=\sum_{\pi\in \Part(k)}f_{\tGa}(\pi) c_\pi(\tmu),$$
 we have $|m_k|\le C^k\tm_k$, and it suffices to prove the lemma for $\tm_k$ instead of $m_k$. We have proved the  claim.
 
 So from now on, we suppose that $\Ga(\cdot)=1$. It allows us to choose a particular model for the random column vector $U_N$: we choose $U_N$ \st its entries $U_N(1), \ld, U_N(N)$ are i.i.d. random variables whose law is the one of $\sqrt{N}B e^{i\Theta}$, where $B$ has law $(1-1/N)\delta_0+(1/N)\delta_1$,  and is independent of $\Theta$, who has uniform law on $[0,2\pi]$. This vector obviously satisfies Hypothesis  \ref{hyp_NfixeUn}  with $\Ga(\cdot)=1$.

Let us define the function $g_\theta(A):=\ff{N}\Tr e^{\theta A}$ on Hermitian matrices. 
To prove the lemma, it suffices to prove that for any $\tta\in \R$, the sequence $\PUNmu(g_\tta)$ stays bounded as $N\to\infty$, which, by definition of the law $\PUNmu$,  will be implied by the fact that $$\limsup_{N\to\infty}\limsup_{n\to\infty}\ff{N}\E[\Tr e^{\f{\tta}{N}\sum_{i=1}^{N\ti n}X_n^i\cdot U_N^iU_N^{i*}}]<\infty,$$ where for each $n$, the $X_n^i$'s are i.i.d. with law $\mu^{*\ff{n}}$, independent of the $U_N^i$'s, some i.i.d. copies of $U_N$.

 The Golden-Thompson inequality (\cf  \cite[Sect. 3.2]{Tao-matrix-book}) states that for $A_1,A_2$ Hermitian matrices, $\Tr e^{A_1+A_2}\le \Tr e^{A_1}e^{A_2}$. In the case where $A_1, \ld, A_k$ are random, independent, \st for all $i$, $\E[e^{A_i}]=m_iI_N$ with  $m_i\ge 0$,   we can deduce, by induction on $k$,  that \be\la{18911.1}\ff{N}\E[\Tr e^{A_1+\cd+A_k}]\le m_1\cd m_k.\ee
 
 Let us now compute $\E[e^{\f{\tta}{N}X_n\cdot U_NU_N^{*}}]$ for $X_n
$ a $\mu^{*\ff{n}}$-distributed variable independent of $U_N$.
Note first that since the \Lvy measure $\si$ of $\mu$ has compact support, by for instance \cite[Lem. 25.6]{MR1739520}, for all $n$,   the Laplace transform of $\mu^{*\ff{n}}$ is defined as an entire function on $\C$ and given by the fonction $e^{\ff{n}\phi(\cdot)}$, where $$
  \phi(\lam)=\lam\gamma+\int\underbrace{(e^{\lam t}-1-\frac{\lam t}{1+t^2} )\frac{1+t^2}{t^2}}_{:=\f{\lam^2}{2}\trm{ for $t=0$}}\si(\ud t).
$$
Since $U_N$ is a column vector,  we have \be\la{29911} e^{\f{\tta}{N}  X_n\cdot U_NU_N^*}=I_N+\f{e^{\f{\tta}{N}  X_n\|U_N\|^2}-1}{\|U_N\|^2}U_NU_N^*,\ee with the convention $(e^t-1)/t=1$ for $t=0$.   
For any diagonal matrix $D$ with diagonal entries on the unit circle, $DU_N$ has the same law as $U_N$,  hence the expectation of the right-hand term of \eqre{29911} is a diagonal matrix. Moreover,  by exchangeability,  for any   $f$,  $$\E[f(\|U_N \|^2)|U_N(k)|^2]=\ff{N}\E[f(\|U_N \|^2)\|U_N\|^2]\qquad\qquad\trm{($1\le k\le N$)},$$so  
 $$\E[e^{\f{\tta}{N} X_n\cdot U_NU_N^{*}}]=I_N+\ff{N}\E[e^{\f{\tta}{N} X_n\|U_N\|^2}-1]I_N=  I_N+\f{n\E[e^{\ff{n}\phi(\tta\|U_N\|^2/N)}-1]}{Nn}I_N.$$ 
 Hence by \eqre{18911.1}, we have $$\ff{N}\E[\Tr e^{\f{\tta}{N}\sum_{i=1}^{N\ti n}X_n^i\cdot U_N^iU_N^{i*}}]\le \lf(1+\f{n\E[e^{\ff{n}\phi(\|U_N\|^2)}-1]}{Nn}\ri)^{Nn}.$$ As $n\to\infty$ with $N$ fixed, the function $n(e^{\ff{n}\phi(\cdot)}-1)$ converges weakly to $\phi(\cdot)$. Moreover, the random variable $\|U_N\|$ takes a finite set of values, hence $n\E[e^{\ff{n}\phi(\|U_N\|^2)}-1]\to \E[\phi(\|U_N\|^2)]$, so that $$\limsup_{n\to\infty}\ff{N}\E[\Tr e^{\f{\tta}{N}\sum_{i=1}^{N\ti n}X_n^i\cdot U_N^iU_N^{i*}}]\le e^{\E[\phi(\tta\|U_N\|^2/N)]}.$$
 To conclude, it suffices to verify that $\E[\phi(\tta\|U_N\|^2/N)]$ stays bounded as $N\to\infty$. It follows from the fact that  the law of $\|U_N\|^2/N$ is the one of $\sum_{i=1}^N B_i$, where the $B_i$'s are i.i.d. with law $(1-1/N)\delta_0+(1/N)\delta_1$.\end{pr}

\subsection{Convergence   of the mean spectral distribution: extension to the case without moments and continuity of the limit}

\begin{propo}\la{571119h24}a) For any $*$-infinitely divisible law $\mu$, the mean spectral distribution $\bar{\Lam}_N^{(\mu)}$ of a $\PUNmu$-distributed random matrix $M$, defined by $$\qquad\qquad\qquad\qquad\int g(t)\bar{\Lam}_N^{(\mu)}(\ud t)=\E\left[\ff{N}\Tr g(M)\right]\qquad\qquad\trm{ (for any function $g$),}$$  converges weakly, as $N\lto\infty$,  to a \pro measure $\Lam_\Gamma(\mu)$ which depends only on $\mu$ and on the function   $\Gamma$ introduced at Equation \eqre{571117h30}.\smallskip

b) The limit measure $\Lam_\Gamma(\mu)$ depends continuously on the pair $(\mu, \Ga)$  for the weak topology. \en{propo}

\begin{pr}a) Let $(\gamma, \si)$ be the \Lvy pair of $\mu$. Note that if $\sigma$ is compactly supported, the result has been established in Proposition \re{571119h07}, so we shall prove the result thanks to Lemma \re{lemCV6312}.

Let $\varepsilon>0$ and $M>1$ such that $\lambda_\infty:=\sigma(\R\bck \left[-M, M\right])<\varepsilon$. Let $\si_\eps$ (resp. $\si_\infty$) the restriction of $\sigma$ to $\left[-M,M\right]$ (resp. $\R\bck \left[-M, M\right]$). Let $\mu_\infty$ be the compound Poisson law with Fourier transform:
$$
\forall \theta\in\R,\ \ \int e^{i\theta t}\ud\mu_\infty(t)=\exp\left(\int_\R\left(e^{i\theta x}-1\right)\frac{1+x^2}{x^2}\si_\infty(dx)\right)
$$
The \Lvy pair of $\mu_\infty$ is $(\ga_\infty,\si_\infty)$ for some $\ga_\infty\in\R$. Let $\ga_\eps=\gamma-\ga_\infty\in\R$, and let $\mu_\eps$ the infinitely divisible law with \Lvy pair $(\ga_\eps,\si_\eps)$. The law  $\mu_\infty$ is a compound Poisson law whose underlying Poisson variable has a parameter less than $2\lam_\infty$, so by construction of the law $\PUNmu$ for $\mu$ a compound Poisson given at Section \re{COCPL5312}, the expectation of the rank of a $\mathbb{P}_{U_N}^{(\mu_\infty)}$-distributed matrix is not larger than $2\eps N$.

 Since $\mu=\mu_\eps*\mu_\infty$, we have $\PUNmu=\mathbb{P}_{U_N}^{(\mu_\eps)}*\mathbb{P}_{U_N}^{(\mu_\infty)}$ by \eqre{Fourier_loi6711}. The mean spectral law of a $\mathbb{P}_{U_N}^{(\mu_\eps)}$-distributed matrix converges for any $\eps$, so the result holds by application of Lemma \re{lemCV6312}.

b) Let $\mu_n$ be a sequence of infinitely divisible \pro measures that converges to a \pro measure $\mu$
 (by  \cite[Lem. 7.8]{MR1739520}, we know that $\mu$ is infinitely divisible). It is known (see \cite[Th. 1, \textsection 19]{GneKol}) that 
if  $(\gamma_n,\sigma_n)$ denotes the \Lvy pair of $\mu_n$ and $(\gamma,\sigma)$ the one of $\mu$, the weak convergence of $\mu_n$ to $\mu$ is equivalent to the fact that \be\la{eq641216h49}\gamma_n\lto\ga\trm{ and }\si_n(g)\lto \si(g)\trm{ for any continuous bounded function $g$}.\ee  
 
 Let us also consider, for each $n\ge 1$, a sequence $(U_N^{(n)})_{N\ge 1}$ of random column vectors which satisfies  Hypothesis  \re{hyp_NfixeUn}. We suppose that there is another sequence $(U_N )_{N\ge 1}$ of random column vectors which satisfies  Hypothesis  \re{hyp_NfixeUn} \st the function $\Ga^{(n)}$ associated to $(U_N^{(n)})_{N\ge 1}$ converges pointwise to the function $\Ga$ associated to the sequence $(U_N)_{N\ge 1}$. 
 
 Let us now prove that the sequence $\La_{\Ga^{(n)}}(\mu_n)$ converges weakly to $\La_\Ga(\mu)$.

 Let us first suppose that there exists $M>0$ such that 
$$
\forall n\geq1,\ \ \sigma_n(\left[-M,M\right]^c)=0.
$$
Then, since the cumulants of $\mu_n$ are more or less the moments of $\si_n$ (see Formula \eqre{571119h57} for more details), we have convergence of the cumulants :  for all $k$, $c_k(\mu_n)\ninf c_k(\mu)$. Hence  by Formula \eqre{eq:0519-2ter}, the moments of $\Lam_{\Ga^{(n)}}(\mu_n)$ tend to the ones of $\Lam_\Ga(\mu)$. It implies the weak convergence, since   $\Lam_\Ga(\mu)$ is determined by its moments, as stated in Proposition \re{571119h07}.

Let us now deal with the general case.  It suffices to prove that the Cauchy transform of  $\Lam_{\Ga^{(n)}}(\mu_n)$ converges pointwise to the one of $\Lam_\Ga(\mu)$. So let us fix $z\in \C^+$ and $\eps>0$.  Define the function $g(t):=\ff{z-t}$. As in the proof of a), one can find $M$ \st $\si(\R\bck[-M,M])\le \eps$ and for all $n$, $\si_n(\R\bck[-M,M])\le \eps$. Hence with some  decompositions $\mu=\mu_\eps*\mu_\infty$ and $\mu_n=\mu_{n,\eps}*\mu_{\infty,\eps}$ as above with $\mu_\eps$ and $\mu_{n, \eps}$ having some \Lvy measures supported on $[-M,M]$ and $\mu_\infty$ and $\mu_{n, \infty}$ some compound Poisson measures, we have \be\la{641216h51}|\La_{\Ga}(\mu)(g)-\La_{\Ga}(\mu_\eps)(g)|\le \f{4\eps}{\Im z},\qquad |\La_{\Ga}(\mu_n)(g)-\La_{\Ga}(\mu_{n,\eps})(g)|\le \f{4\eps}{\Im z}.\ee
This is a consequence of the resolvent identity (see Equation (\ref{residentity_6412}), Lemma \ref{lemCV6312} in the Appendix).
Moreover, by what \eqre{eq641216h49}, $\mu_{n,\eps}$ converges weakly to $\mu_\eps$ as $n\to\infty$. So by what precedes about the compact case,  $\La_{\Ga^{(n)}}(\mu_{n,\eps})$ converges weakly to $\La_{\Ga }(\mu_{ \eps})$. Joining this to \eqre{641216h51}, we get the desired result. 
\end{pr}

\subsection{Concentration of measure and almost sure convergence}
To conclude the proof  of  Theorem \re{mainTh5711}, we still need to prove that that not only the {\emph{mean spectral law}} of a $\PUNmu$-distributed random matrix $M$ converges to $\Lam_\Ga(\mu)$ as $N\lto\infty$, but also, almost surely, the {\emph{empirical spectral law}}. By the Borel-Cantelli lemma, it follows directly from what precedes and  the following lemma. 

 \begin{rmk}[The associated \Lvy process]\la{LevyPro7711}{\rm   
 There exists a \Lvy process $(M_s)_{s\geq0}$ with values in the space of $N\ti N$ Hermitian matrices such that $M_s$ has distribution $\mathbb{P}_{U_N}^{(\mu^{*s})}\equiv({\PUNmu})^{*s}$ for every $s\geq 0$ (this is due to the fact that the distribution $\PUNmu$ is infinitely divisible, see \eg   \cite[Cor. 11.6]{MR1739520}). Let $\left(P_s^{(N,\mu)},s\geq 0\right)$ be the corresponding semi-group: for every Hermitian matrix $A$, every bounded function $f$:
$$
P_s^{(N,\mu)}(f)(A)=\E\left[f\left(M_s+A\right)\right].
$$
Let $\A^{(N,\mu)}$ the infinitesimal generator associated to this
semi-group. Its domain contains the set of  twice continuously
differentiable functions  on $N\ti N$ Hermitian matrices vanishing at
infinity and if $(\gamma, \si)$ denotes the \Lvy pair of $\mu$ (\emph{cf.} Proposition-Definition \ref{definfdivclassique}), for any Hermitian matrix $A$,  we have, by \cite[Th. 31.5]{MR1739520},
\beq
\A^{(N,\mu)}(f)(A)&=&\gamma \ud f(A)\left[D_N\right]\\ && +\int\left\{N\left[\E\left[f\left(A+xU_NU_N^*\right)\right]-f(A)\right]-\frac{x}{1+x^2}\ud f(A)\left[D_N\right]\right\}\frac{1+x^2}{x^2}\sigma(\ud x),
\eeq where $\ud f$ denotes the derivative of $f$ and $D_N=\mathbb{E}\left[U_NU_N^*\right]$. 
  Note that for $f: \R\to \R$  extended   to Hermitian matrices  by  spectral calculus, we have  $\ud f(A)\left[I_N\right]=f'(A)$.}\en{rmk}

\begin{lem}\label{lem:0519-1} 
Consider $\mu$ an infinitely divisible law and   $g:\R\to\R$ a Lipschitz function with finite total variation. Then for every $\varepsilon>0$, there exists $\delta>0$ such that for all $N\geq 1$,
$$
\Pro\left(\left\vert\ff{N}\Tr g(M)-\E\left[\ff{N}\Tr g(M)\right]\right\vert>\varepsilon\right)\leq 2e^{-N\delta}
$$
with $M$ a $\PUNmu$-distributed random matrix.
\end{lem}

\beg{pr} This is an extension of Theorem III.4 in \cite{cab-duv-BP} which is established for a special case of $U_N$. Its proof only uses the fact that $\Vert U_N/\sqrt{N}\Vert_2=1$ (one has to notice that in the present paper, what  plays the role of $U_N$ in \cite{cab-duv-BP} is $U_N/\sqrt{N}$). Theorem III.4 in \cite{cab-duv-BP} can be readily extended to the case when $\E\left[\Vert U_N/\sqrt{N}\Vert_2^4\right]$ is bounded w.r.t. $N$. Then, we conclude noticing that by exchangeability and Equation \eqre{571117h30}, $\E\left[\Vert U_N/\sqrt{N}\Vert_2^4\right]$ converges to $\Gamma(2)+\Gamma(1,1)$. Indeed, \beq \E\left[\Vert U_N/\sqrt{N}\Vert_2^4\right]&=&N^{-2}
\sum_i\E[|U_N(i)|^4]+N^{-2}\sum_{i\ne j}\E[|U_N(i)|^2|U_{N}(j)|^2]\\
&=& N^{-1}\E[|U_N(1)|^4]+(1-N^{-1})\E[|U_N(1)|^2|U_{N}(2)|^2]
\eeq
\en{pr}

\section{Proof of Proposition \re{unbounded_support_121011}}\la{pr_unbounded_support_121011}
For each $n\ge 1$, let $m_n$ be the $n^{th}$ moment of $\Lam_\Ga(\mu)$.  Let $c=\inf \Gamma(n)^{1/n}$. It suffices to prove that there is $\eps>0$ \st for all $n$ even, $$m_n\ge \eps^n\ti\trm{$n^{th}$ moment of a standard Gaussian variable},$$ \ie that for all $n$ even, $m_n\ge \eps^n|\Part_2(n)|,$ where $\Part_2(n)$ is the number  of pairings  of $\{1, \ld, n\}$. The formula of $m_n$ is given by \eqre{eq:0519-2ter}, where each term is positive. Moreover, we know that $c_2(\mu)=\Var(\mu)>0$. Hence it suffices to notice that for any $k\ge 1$ and any $\pi\in \Part_2(k)$, $f_\Ga(\pi)\ge c^k$, which follows from the expression of $f_\Ga(\pi)$ as a sum of positive terms  in \eqre{20911.19h}, where the term associated to the trivial partition $\tau$ with one block is $$\prod_{J\in \pi}\Ga(|J|)\ge \prod_{J\in \pi}c^{|J|}=c^k.$$

 \section{Proof of Theorem \re{HT5412}}
 First, in the case where $\al>1$, we can suppose the $Y_{ij}$'s to be centered. Indeed, replacing  the $Y_{ij}$'s by $Y_{ij}-\E[Y_{ij}]$ is  a rank-two perturbation of $M_{N,p}$ and by Lemma \re{lemCV6312}, a rank-two perturbation of an Hermitian matrix has no influence on the weak convergence of its spectral law.
 
Let  $(P(t ), t\ge 0) $ be a standard    Poisson process, independent of the other variables. For each  $B>0$, let us introduce both following approximations of the random matrix $M_{N,p}$ :     $$\widehat{M}_{N,p}:=\ff{a_N^2}\sum_{j=1}^{p} X_j\cdot \widehat{V}^j(\widehat{V}^j)^*\qquad\trm{ and }\qquad\widetilde{M}_{N}:=\ff{a_N^2}\sum_{j=1}^{P(N\lam)} X_i\cdot \widehat{V}^{j}(\widehat{V}^j)^*$$ with $\widehat{V}^j:=(Y_{ij}\one_{|Y_{ij}|\le Ba_N})_{j=1}^N \in \K^{N\ti 1}$ (this column vector depends implicitly on $N$ and on the cutoff parameter $B$). It is noticed in Sections 1 and 9 of \cite{BAGheavytails} that for $Z^{(B)}_N:=\f{\sqrt{N}}{a_N}Y_{11}\one_{|Y_{11}|\le Ba_N}$, we have, for each $n\ge 1$, \bes\la{explosion_moments_121011} \f{\E[|Z_N^{(B)}|^{2n}]}{N^{{n}-1}}\Ninf \f{\al}{2n-\al} B^{2n-\al} \ees and $\sqrt{N}\E[Z_N^{(B)}]=O(1)$ (this is where the recentering is necessary when $\al>1$).     So Hypothesis     \re{hyp_NfixeUn}   holds for $U_N$ distributed as the $\widehat{V}^j$'s, with \be\la{GaB6412}\Ga(n_1, \ld, n_k)=\Ga^{(B)}(n_1, \ld, n_k):= B^{2(n_1+\cdd+n_k)-k\al}\prod_{\ell=1}^k\f{\al}{2n_\ell-\al}.\ee
 As a consequence, by Theorem \re{mainTh5711}, the empirical spectral law of $\widetilde{M}_{N}$ converges almost surely, as $N\to\infty$, to the law $\Lam_{\Ga^{(B)}}(\mu)$, with $\mu$ the (compound Poisson) law of $\sum_{j=1}^{P(\lam)}X_j$. By Lemma \re{lemCV6312} (applied with $M_N=\widehat{M}_{N,p}$ and and $M_N^\eps=\widetilde{M}_{N}$ for all $\eps>0$), the same holds for $\widehat{M}_{N,p}$, because 
 by the Law of Large Numbers, as $N,p\to\infty$, with $p/N\to\lam$, $$\qquad\qquad\f{\op{rank}(\widetilde{M}_{N}-\widehat{M}_{N,p})}{N}\lto 0 \qquad\trm{(almost surely)}.$$ 
 
 To prove the convergence of the empirical spectral law of $M_{N,p}$, by Lemma \re{lemCV6312}, it suffices to prove that for any $\eps>0$, if $B$ is large enough, then almost surely, for $N,p$ large enough, \be\la{6412.12h44} \op{rank}(M_{N,p}-\widehat{M}_{N,p}) \le N\eps.\ee But we have $$\op{rank}(M_{N,p}-\widehat{M}_{N,p})\le \sum_{j=1}^p 2\one_{V^j\ne\widehat{V}^j}.$$The above right-hand-term is a sum of $p$ independent Bernoulli variables with parameter $$q_N:=1-(\mathbb{P}(|Y_{11}|\le Ba_N))^N.$$We have, for $L$ the slow variations function introduced at \eqre{queue_distrib_6412},  $$\mathbb{P}(|Y_{11}|\le Ba_N)=1-B^{-\al}a_N^{-\al}L(Ba_N)=1-B^{-\al}a_N^{-\al}L(a_N)\ti\f{L(Ba_N)}{L(a_N)}=1-\f{B^{-\al}(1+\epsilon_N)}{N},$$ where $\epsilon_N$ is a sequence tending to zero. As a consequence, as $N\to\infty$, $q_N$ tends to $1-e^{-B^{-\al}}$. By some Laplace transform estimates like p. 722 of \cite{BAGheavytails}, there is  a constant $c$ such that $$\mathbb{P}(\op{rank}(M_{N,p}-\widehat{M}_{N,p})\ge 4pq_N)\le e^{-cpq_N},$$ hence by Borel-Cantelli's lemma, almost surely, for $N,p$ large enough ($N$ and $p$ grow together in such a way that $p/N\to\lam$), $\op{rank}(M_{N,p}-\widehat{M}_{N,p})\le 4pq_N$. But $1-e^{-B^{-\al}}\sim B^{-\al}$ as $B\to\infty$, so for any fixed $\eps$, one can choose $B$ \st almost surely, for $N,p$ large enough, \eqre{6412.12h44} holds. It concludes the proof of the convergence. 
 
 Let us now denote by $\nu$ the law of the $X_i$'s and by $\La_\al(\nu, \lam)$ the limit of the empirical spectral law of $M_{N,p}$. We shall prove that this is a continuous function of $(\nu, \al, \lam)$.  By the formula given for $\Ga^{(B)}$ at \eqre{GaB6412}, for any fixed cutoff $B$, the limit spectral law of $\widehat{M}_{N,p}$ depends continuously on the parameters $(\nu, \al, \lam)$. Since the cutoff necessary to obtain a right-hand-term $N\eps$ in \eqre{6412.12h44} can be chosen uniformly on $\al$ and $\lam$ as soon as they vary in compact sets which do not contain zero (regardless to the law of the $X_j$'s), this proves the continuity.

\section{Proof of Proposition \re{871119h30}}\la{871119h31}

We first treat (together)  the Gaussian and uniform cases.  Exchangeability is obvious. 
Let us now prove 
that for any diagonal matrix $D=\diag(\eps_1,\ld, \eps_N)$ with diagonal entries in $\{z\in \K\ste |z|=1\}$, \be\la{L111011}U_{N}U_{N}^*\eqlaw DU_{N}U_{N}^*D^*.\ee It will imply that the expectation of \eqre{571117h30_mjj} is null   as soon as the condition on multisets given at Remark \re{tildeU_N} is satisfied.  
Let $\{e_1, \ld, e_N\}$ denote the canonical basis of $\C^{N\ti 1}$. For each $\ell=1, \ld, N$,  let us denote by $\mc{L}_\ell$ the law of $U_N$ conditionally to the event $X_N=e_\ell$ ($X_N$ is the random variable with uniform law on $\{e_1, \ld, e_N\}$ that was introduced after \eqre{defGauss_111011_71012}). We know that the law $\mc{L}$ of $U_N$ is given by  \be\la{169111}\mc{L}(B)=\ff{N}\sum_{\ell=1}^N \mc{L}_\ell(B)\ee for any Borel set $B\subset \K^N$. Hence it suffices to prove that for all $\ell$, an  $\mc{L}_\ell$-distributed vector $V$ satisfies 
 \eqre{L111011}. So let us fix $\ell\in \{1, \ld, N\}$, consider such a vector $V=(V(1), \ld, V(N))^T$ and such a matrix  $D=\diag( \eps_1,\ld, \eps_N)$. Note that the entries of $V$ are independent and that all of them except the $\ell$-th one are invariant, in law, by multiplication by any $\eps_i$. Hence for $D'=\eps_\ell^{-1}I_N$, we have $$(D'D)V\eqlaw V.$$ As a consequence, we have $$(D'D)VV^*(D'D)^*\eqlaw VV^*.$$ But we also have $(D'D)VV^*(D'D)^*=DVV^*D^*$, hence \eqre{L111011} is proved. 
It remains to prove \eqre{571117h30} and more precisely Formulas (\ref{871119h29gauss}), (\ref{871119h292gauss}) and (\ref{871119h293gauss}). So let us fix $k\ge 1$, $n_1, \ld, n_k\ge 1$ and find out the limit, as $N\to\infty$,  of $$N^k\E[|U_N(1)|^{2n_1}\cd |U_N(k)|^{2n_k}].$$ Note that by \eqre{169111}, we have \beqy\la{169112}&N^{k-(n_1+\cd+n_k)}\E[|U_N(1)|^{2n_1}\cd |U_N(k)|^{2n_k}]=&\\ 
\la{1691122} & \lf(1-\f{k}{N}\ri)N^{k-(n_1+\cd +n_k)}(1-e^{-t})^{n_1+\cd +n_k}\prod_{j=1}^k\E[ |g_1|^{2n_j}]
 +&\\ \la{1691123} & \sum_{j_0=1}^kN^{k-(n_1+\cd +n_k)+n_{j_0}-1}(1-e^{-t})^{n_1+\cd +n_k-n_{j_0}}\E[|e^{-\f{t}{2}}+\f{g_1}{\sqrt{N}}|^{2n_{j_0}}]\prod_{\substack{j=1\\ j\ne j_0}}^k\E[ |g_1|^{2n_j}].&
\eeqy
First, as soon as one of the $n_j$'s is $\ge 2$, the  term of   \eqre{1691122} vanishes as $N\to\infty$. The same happens for each of the terms of the sum \eqre{1691123} if two of the $n_j$'s are $\ge 2$. Equation \eqre{571117h30}  follows easily. 

Let us now consider the case where $U_N=U_{t}$, for $(U_t/\sqrt{N})_{t\ge 0}$ a solution of \eqre{EDSbro} whose initial law is the uniform law on the canonical basis.    Let us introduce the $N\ti N$   matrix $P_{t}$ (the dependence on
$N$ is implicit in $P_t$ and in $U_t$) 
$$
P_t:=U_tU_t^*/N.
$$
By a direct application of the matricial It\^o calculus (see \cite[Sect. 2.1]{flobulletin}), the process $(P_t)_{t\ge 0}$ satisfies the SDE 
\be\la{SDEP_t}
\ud P_t=(\ud K_t)P_t-P_t(\ud K_t) +\left(\ff{N}I-P_t\right)\ud s.
\ee First, from this SDE, it follows easily that $\Tr P_t=1$ for all $t$, so that $\|U_t\|^2=N$. Second, 
it is easy to see that the law of $K$ is invariant under conjugation by any unitary matrix. By uniqueness, in law, of the solutions  of \eqre{EDSbro} and \eqre{SDEP_t}, given the initial conditions, it follows that for any permutation matrix $Q$ and any matrix $D$ as in \eqre{L111011},  $QU_t\eqlaw U_t$ and $DU_tD^*\eqlaw P_t$.  It remains to prove \eqre{571117h30} and Formulas (\ref{871119h29gauss}), (\ref{871119h292gauss}) and (\ref{871119h293gauss}). 
 For each $n_1, \ld, n_k\ge 0$, let us define 
$$
\tGt(n_1,\ldots,n_k):=N(N-1)\cdots(N-p+1)\E\left[|U_t(1)|^{2n_1}\cd\cd|U_t(k)|^{2n_k}\right]N^{-(n_1+\cd+n_k)},
$$ 
where $p={\Card{\{\ell\ste n_\ell\ne 0\}}}$ (the dependence of  $\tGt$ on $N$
is implicit). We have to prove that  
$$
\tGt(n_1,\ldots,n_k)\Ninf \Gamma_t(n_1,\ldots,n_k).
$$ 
Due to exchangeability and to $\|U_t\|=N$, we have $\tGt(1)=1$ and  
\be\la{871119h293dafter2}
\tGt(n_1,\ldots,n_k)=\tGt(n_1+1,n_2,\ldots,n_k)+\cdots+\tGt(n_1,\ldots,n_{k-1},n_k+1)+\tGt(n_1,\ldots,n_k,1),
\ee
 and since all $|U_t(i)|$'s are $\le 1$, if one considers two families
 $n_1, \ld, n_k\ge 1$ and  $n_1', \ld, n_k'\ge 1$ \st $n_j\le n_j'$ for all
 $j$, then 
\be\la{871119h293dafter} 
\tGt(n_1,\ldots,n_k)\ge \tGt(n_1',\ldots,n_k').\ee 
It is easy to see that thanks to  
  \eqre{871119h293dafter2} and \eqre{871119h293dafter} and to the obvious formula $\tGt(1)=1$,  it suffices to prove that 
\be\la{97111}
\tGt(2,2)\Ninf 0 \quad\quad\trm{ and }\quad\quad  \tGt(n)\Ninf e^{-nt} \qquad \trm{ for all
  $n\ge 2$.}
\ee
For any multiset $\{k_1, \ld, k_n\}$ of integers in $\{1, \ld, N\}$, let us define $$\ga_t(k_1, \ld, k_n)=N^{-n}|U_{t}(k_1)|^2\cdots |U_{t}(k_n)|^2 $$ (the dependence of $\ga_t(\cd)$ in $N$ is implicit). 
By exchangeability, $\E\left[\ga_t(k_1, \ld, k_n)\right]$ only
depends on the level sets partition of the map $j\mapsto k_j$. To prove
\eqre{97111}, we have to prove  
\be\la{97111bis}
N^2\E\left[\ga_t(1,1,2,2)\right]\Ninf 0
\ee 
and   
\be\la{97117h24bis}
\qquad\qquad\qquad\qquad N\E\left[\ga_t(\underbrace{1,\ld, 1}_{n\trm{ times}})\right]\Ninf
e^{-nt} \qquad\qquad \trm{ for all $n\ge 2$.}
\ee
Recall that $\{e_1, \ld, e_N\}$ denotes the canonical basis of $\C^{N\ti 1}$. For each $k\ge 1$,  
$$
\ud \ga_t(k)=e_k^*(\ud K_t)P_te_k-e_k^*P_t\ud K_te_k +\left(\ff{N}-\ga_t(k)\right)\ud
t.
$$ 
As an application  of the matricial It\^o calculus (see \cite[Sect. 2.1]{flobulletin}),  we get the quadratic variation dynamics
$$
\ud\lan \ga_t(k), \ga_t(\ell)\ran=
-\f{2}{N}\left[\ga_t(k)\ga_t(\ell)-\one_{k=\ell} \ga_t(k)\right]\ud t.
$$ 
By It\^o's formula for products of semi martingales, it follows  that 
$\E\left[\ga_t(k_1, \ld, k_n)\right]$ is a smooth
function of $t$ which satisfies the differential system 
\beqy\la{971108h27} 
\qquad\partial_t\E\left[ \ga_t(k_1, \ld, k_n)\right] &=&
-n\lf(1+\f{n-1}{N}\ri)\E\left[\ga_t(k_1, \ld, k_n)\right]\\ 
\nonumber&&+\ff{N}\sum_{j=1}^n(2\Card{\{i\ste i<j, k_i=k_j\}}+1)\E\left[\ga_t(k_1,\ld,\hk_j,\ld, k_n)\right],
\eeqy
where we used the convention $\ga_t(k_1, \ld, k_n)=1$ if $n=0$.
Let us now prove \eqre{97117h24bis}. For each $n\ge 1$, we define 
$$
q_t(n)=Ne^{nt}\E\left[\ga_t(\underbrace{1, \ld, 1}_{n\trm{ times}})\right]
$$ 
(the dependence of $q_t(n)$ on $N$ is implicit). 
Note first that since $\|U_t\|^2=N$,  for all $n\ge 1$,   
$$
\E\left[|U_{t}(1)|^{2n}\right]N^{-n+1}\le  \E\left[|U_{t}(1)|^2\right]=1,
$$ 
so  that  $q_t(n)\le e^{nt}$. Moreover, by the differential system of
\eqre{971108h27}, we have, for all $n\ge 2$, 
$$
\partial_t q_t(n)=-\f{n-1}{N}q_t(n)+\ff{N}\sum_{i=1}^n(2j-3)q_t(n-1).
$$ 
It follows that $\partial_t q_t(n)$ converges to zero uniformly on every
compact subset of $\R_+$ as $N\lto\infty$. Since for $t=0$, $q_t(n)=1$ for
each $n$, we get that for each $n\ge 2$,  
$$
q_t(n)\Ninf 1,
$$ 
so that \eqre{97117h24bis} is proved.
Let us at last prove \eqre{97111bis}. Thanks to the differential system of
\eqre{971108h27} and to the fact that
$\E\left[\ga_t(1,1,2)\right]=\E\left[\ga_t(1,2,2)\right]$, we have 
$$
\partial_t\E\left[\gamma_t(1,1,2,2)\right]=-4\left(1+\frac{3}{N}\right)\E\left[\gamma_t(1,1,2,2)\right]+\frac{8}{N}\E\left[\gamma_t(1,1,2)\right]
$$
Due to exchangeability and to  $\|U_t\|=N$,  both $\E\left[\gamma_t(1,1,2)\right] $ and $  \E\left[\gamma_t(1,1,2)\right]$ are not greater than $\frac{1}{N(N-1)}$. Therefore, we get:
$$
\left\vert\partial_t\left(N^2e^{4t}\E\left[\gamma_t(1,1,2,2)\right]\right)\right\vert\leq\frac{11}{N}
$$
Since $\gamma_0(1,1,2,2)=0$, this proves \eqref{97111bis}.\qed

\section{Proof of Proposition \ref{propo:kappa}}
\label{sec:kappa}

Let $\pi\in\Part(k)$ and $\Gamma:=\Gamma_t$; we have to prove that $f_{\Ga_t}(\pi)=e^{\kappa(\pi)t}$. By  the very definition of $\kappa$, it is  
obvious that $\kappa(\pi)=\kappa(\thin(\pi))$ and that $\kappa$ is additive on the connected components of $\pi$. 
Therefore, due to Lemma \ref{pte:0213-3gusgus}, it is enough to
prove that $f_{\Ga_t}(\pi)=e^{\kappa(\pi)t}$ for $\pi$ thin and connected. Note
that in this case, $\kappa(\pi)$ is equal to $k$. 
 
We start from Formula \eqre{20911.19h} for $f_{\Ga_t}(\pi)$ and  are going to prove that the only non-vanishing term in this formula
corresponds to the trivial partition $\tau$ with only one block.  This will imply:
$$
f_{\Ga_t}(\pi)=\prod_{J\in\pi}\Gamma_t(\vert J\vert)=\prod_{J\in\pi}e^{-\vert
  J\vert t}=e^{-kt}=e^{-\kappa(\pi)t}
$$
and this will prove the proposition.

For $k=1$, it is trivial. 
Let us suppose $k>1$. Since $\pi$ is thin and connected, it has a crossing
(and this implies in fact $k\geq4$). 
Let us consider a non trivial $\tau$ (i.e. with at least two blocks) \st $H:=H(\pi, \tau)$ has no cycle. Let $c$ be a coloring of the edges of $\Gspi$ with kernel $\tau$.

Since $H$ is connected, there
  exists two blocks $E$ and $E'$ of $H$ with a non-empty intersection.     Let
  $J_0\in\pi$ such that $\{J_0\}=E\cap E'$ (since $H$ has no cycle, this
  intersection is a singleton). Since $\pi$ is thin, by Proposition \re{prop_fond_real_14911}, $E$ is not a
  singleton. Let $J_0'\in E\setminus\{J_0\}$. Since  the graph $G_{/\pi}$ is a circuit, there exist at least one path from $J_0$
to $J_0'$ and one path from $J_0'$ to $J_0$. But if they were the only ones,   there would be no crossing between $J_0$ and the others blocks of
$E$, and $\pi$ would not be connected (\emph{cf.} for instance the
(counter-)example page \pageref{fig:vraigraphequotient}). Therefore there exist at least two paths from $J_0$
to $J_0'$ and two paths from $J_0'$ to $J_0$. Hence we deduce that $n_c^+(J_0,i)\geq
2$ for a certain colour $i$ (see page \pageref{ncdefplus} for the definition of $n_c^+(J_0,i)$). With the same argument, we get $n_c^+(J_0,j)\geq2$ for a certain $j\ne i$. This implies that 
$$
\Gamma(n_\tau(J_0))=\Gamma(n_c^+(J_0,\ell),\ell\ge 1)=0.
$$
Hence, $\prod_{J\in \pi}\Gamma(n_\tau(J))=0$, and this implies that the only $\tau$ with a non-vanishing contribution to the
computation of $f_{\Ga_t}(\pi)$ is the trivial one with only one block. \qed

\section{Proof of Proposition \re{771113h50}}\la{771113h51}

\beg{lem}For any Hermitian matrices $A,B$ with $B$ positive, any   column vector $V$, any constant $c\in\R$:
\begin{equation}\label{eq:0519-4}
\E\left[(A+iB-cX_t VV^*)^{-1}\right]=(A+iB+it\vert c\vert VV^*)^{-1}
\end{equation}
with $X_t$ a random variable with law $\Cc_t$.\en{lem}

\beg{pr} On can suppose $V$ to have unit norm.
First, for $f$   a rational function bounded with non positive degree  whose poles are all in $\C^+:=\{z\in\C\ste \Im(z)>0\}$,   $\Cc_t(f)=f(-t i)$ (it can be proved via the residue formula if $f(z)=\ff{z-a}$ and can then be generalized by density). Second, for $Q$   the orthogonal projector on $V^\bot$,
$$
\mathrm{det}\left(A+iB-cx VV^*\right)=(\langle V,(A+iB)V\rangle-cx)\mathrm{det}_{\vert Q}(Q(A+iB)Q)
$$
with $\mathrm{det}_{\vert Q}$ the determinant on the space $V^\bot$. If
we suppose $c>0$, this implies that the entries of $(A+iB-cx VV^*)^{-1}$
are rational functions of $x$ with poles in $\C^+$. The equality
(\ref{eq:0519-4}) can then be deduced from these two remarks.
\en{pr}

\noindent{\bf Proof of Proposition \re{771113h50}.} Let $X_{n}^1,\ldots,X_{n}^{Nn}$, $U_{N}^{1},U_N^2, \ldots$ be 
as in Theorem \ref{defPUNmu}, with $X_{n}^{i}\sim\Cc_{\frac{t}{n}}$.   If we apply repeatedly the preceding equality, we obtain that for any Hermitian matrix $A$, $z\in\C^+$, $c\in\R$:
$$
\E\left[(zI_N-A- (X_{n}^1U_{N}^{1}U_{N}^{1*}+\cdots+X_{n}^{Nn}U_{N}^{Nn}U_{N}^{Nn*})/N)^{-1}\vert U_{N}^{1},\ld,U_N^{Nn}\right]$$ $$
=\left(zI_N+it\frac{ U_{N}^{1}U_{N}^{1*}+\cdots+ U_{N}^{Nn}U_{N}^{Nn*}}{nN}-A\right)^{-1}
$$
When $n\lto\infty$, $\frac{ U_{N}^{1}U_{N}^{1*}+\cdots+ U_{N}^{Nn}U_{N}^{Nn*}}{nN}$ tends a.s. to $\E\left[U_NU_N^*\right]=I_N$ by the law of large numbers, and $(X_{n}^1U_{N}^{1}U_{N}^{1*}+\cdots+X_{n}^{Nn}U_{N}^{Nn}U_{N}^{Nn*})/N$ tends in law to $\mathbb{P}_{U_N}^{(\Cc_t)}$ following Theorem \ref{defPUNmu}. Therefore, it is now easy to deduce that:
$$
\E\left[(zI_N-A-M_t)^{-1}\right]=(zI_N+it I_N-A)^{-1}=P_{t}(f_z)(A)
$$
with $f_z(x)=(z-x)^{-1}$. This establishes the proposition for $f=f_z$. Using routine density and linearity arguments, this is enough to prove the proposition.\qed

\section{Appendix}

\subsection{Matrix approximations in the sense of the rank}  
The \emph{Cauchy transform} of  a finite measure  $\mu$ on $\R$ is $G_\mu(z):=\int \f{\mu(\ud t)}{z-t}$, for  $z\in\C^+$. 
\beg{lem}\la{lemcauchy6412}Let $\mu_N$ be a sequence of random \pro measures on the real line \st for any $\eps>0$, there is another sequence $\mu_N^\eps$   converging weakly to a deterministic \pro measure $\mu^\eps$ and \st almost surely, for $N$ large enough,  for any $z\in \C^+$,  \be\la{63121}\lf|G_{\mu_N}(z)-G_{\mu_N^\eps}(z)\ri|\le \f{\eps}{\Im z}.\ee Then $\mu_N$ converges weakly to a deterministic \pro measure $\mu=\lim_{\eps\to 0}\mu^\eps$. 
\en{lem}

\beg{pr} Almost surely, by \eqre{63121}, for any $z\in \C^+$, the sequence $G_{\mu_N}(z)$ is a Cauchy sequence, hence converges to a deterministic limit $G(z)$ \st for any $\eps$, \be\la{63121.bis}\lf|G_{\mu}(z)-G_{\mu^\eps}(z)\ri|\le \f{\eps}{\Im z}.\ee So by \cite[Th. 2.4.4 b)]{agz09},  $\mu_N$ converges vaguely  to a   measure $\mu$ with total mass $\le 1$ \st $G_\mu=G$ and $\mu=\lim_{\eps\to 0}\mu^\eps$.    To see that $\mu$ is a \pro measure, it suffices to notice that for any      measure $\nu$,   $iy\Im(G_\nu(iy))$ increases to $\nu(\R)$ as $y\uparrow\infty$, and then to use \eqre{63121.bis}.
\en{pr}

For $M$ a random $N\ti N$ Hermitian matrix with eigenvalues $\lam_1, \ld, \lam_N$, recall that the \emph{empirical spectral law} of $M$ is the random \pro measure $\mu_M$ defined by \be\la{710129h45}\mu_M:=\ff{N}\sum_{i=1}^N\delta_{\lam_i}\ee and that the  \emph{mean spectral law} of $M$ is the deterministic  \pro measure $\ovl{\mu}_M$ defined, for any bounded Borel function $f$, by  \be\la{710129h452}\ovl{\mu}_M(f):=\E(\mu_M(f)],\ee
with $\mu_M$ as in \eqre{710129h45}.

\beg{lem}\la{lemCV6312}For each $N$, let $M_N$ be an $N\ti N$ random Hermitian matrix. Suppose that for any $\eps>0$, there is a sequence $M_N^\eps$ of random Hermitian matrices whose empirical (resp. mean) spectral law converges almost surely (resp. converges) to a deterministic \pro measure $\mu^\eps$ and \st almost surely, for $N$ large enough, $\op{rank}(M_N-M_N^\eps)\le N\eps$ (resp. \st for $N$ large enough,   $\E[\op{rank}(M_N-M_N^\eps)]\le N\eps$). Then the empirical (resp. mean) spectral law of $M_N$ converges almost surely (resp. converges) to a deterministic \pro measure $\mu=\lim_{\eps\to 0}\mu^\eps$
\en{lem}

\beg{pr} Let us first prove the \emph{almost sure version} of the lemma. For any Hermitian matrix $H$ and any $z\in \C^+$,      we have  $\|(z-H)^{-1}\|_{\op{op}}\le (\Im z)^{-1}$. So by the formula $$(z-B)^{-1}-(z-A)^{-1}=(z-B)^{-1}(B-A)(z-A)^{-1},$$ we have  \begin{equation}\label{residentity_6412}
\lf| \Tr ({z-M_N})^{-1}- \Tr ({z-M_N^\eps})^{-1}\ri|\le {2}({\Im z})^{-1} \ti \op{rank}(M_N-M_N^\eps)
\end{equation}   and for any $\eps$, with \pro one, for $N$ large enough, for all $z\in \C^+$ $$
\lf| \ff{N}\Tr ({z-M_N})^{-1}- \ff{N}\Tr ({z-M_N^\eps})^{-1}\ri|\le \f{2\eps}{\Im z},\quad\trm{ \ie }\quad |G_{\mu_{M_N}}(z)-G_{\mu_{M_N^\eps}}(z)|\le \f{2\eps}{\Im z}.$$
It follows that the sequences $\mu_N:=\mu_{M_N}$ and $\mu_N^\eps:=\mu_{M_N^\eps}$ satisfy the hypotheses of the previous lemma, which allows to conclude.

Let us now prove the \emph{expectations version} of the lemma. Taking the expectation of \eqre{residentity_6412}, we get $$|G_{\ovl{\mu}_{M_N}}(z)-G_{\ovl{\mu}_{M_N^\eps}}(z)|\le \f{2\eps}{\Im z}.$$
It follows that the (non random) sequences $\mu_N:=\ovl{\mu}_{M_N}$ and $\mu_N^\eps:=\ovl{\mu}_{M_N^\eps}$ satisfy the hypotheses of the previous lemma, which allows to conclude.
\en{pr}

\subsection{Infinitely divisible distributions and cumulants}\la{appendiceinfdiv}
\subsubsection{The classical case}
Let us   recall the definition of infinitely divisible laws with respect to the classical convolution $*$ (see \cite{petrov}).
\begin{propdef}\la{definfdivclassique}Let $\mu$ be a \pro measure on $\R$. The we have equivalence between : 
\bgt\ite[(i)] there is a sequence $(k_n)$ of positive integers tending to $+\infty$ and a sequence $(\nu_n)$ of \pro measures on the real line \st as $n\lto \infty$,   \be\la{5611_20h48}\underbrace{\nu_n*\cdots*\nu_n}_{\trm{$k_n$ times}}\lto \mu,\ee
\ite[(ii)] there is a continuous  semigroup  $(\mu^{*t})_{t\in [0, \infty)}$  for the convolution  $*$ starting at $\delta_0$  \st $\mu^{*1}=\mu$,
\ite[(iii)]   the Fourier transform  of $\mu$ has the form $\int_{t\in \R} e^{it\xi}\ud\mu(t)=e^{\Psi_\mu(\xi)}$, with \be\la{66117h29}\Psi_\mu(\xi)= i\gamma \xi+\int_\R\underbrace{(e^{it\xi}-1-\f{it\xi}{t^2+1})\f{t^2+1}{t^2}}_{:=-\f{\xi^2}{2}\trm{ for $t=0$}}\ud \sigma(t),\ee  $\gamma$ being a real number and $\si$ being a finite positive measure on $\R$.
\ent
\end{propdef}
In this case, $\mu$ is said to be \emph{$*$-infinitely divisible}, the pair $(\gamma, \si)$, called the \emph{L\'evy pair} of $\mu$,  is unique, and $\mu$ will be denoted by $\nu_*^{\gamma,\si}$.\\

Part (i) of the above definition characterizes such laws as the  limit laws of sums of i.i.d. random variables, Part (ii) expresses them as the distributions of one-dimensional marginals of L\'evy processes, and Part (iii) is known as the \emph{L\'evy-Kinchine formula}. The pair $(\gamma, \si)$  can be interpreted as follows: $\gamma$ is a drift, $\si(\{0\})$ represents the brownian component of the L\'evy process 
 associated to $\mu$, and the measure $\one_{x\ne 0}\f{1+x^2}{x^2}\ud \si(x)$, when it is finite, represents its Poisson compound part (when this measure is not finite, the \Lvy process can be understood as a limit of such decompositions).     
 
 \subsubsection{The free case}
 In \cite{BV92,defconv}, Bercovici and Voiculescu have proved that Proposition-Definition \re{definfdivclassique} stays true if one replaces the classical convolution $*$ by the free additive convolution $\bxp$ and Formula \eqre{66117h29} by the following formula for the $R$-transform of $\mu$:  \be\la{66117h292} R_\mu(z)=\gamma +\int_\R\f{z+t}{1-tz}\ud \sigma(t)\ee (recall that the \emph{$R$-transform} of $\mu$ is defined by the formula $R_\mu(z)=G_\mu^{-1}(z)+\ff{z}$, see \cite{ns06,agz09}). In this case,   $\mu$ is said to be \emph{$\bxp$-infinitely divisible} and is  denoted by $\nu_\bxp^{\gamma, \si}$.\\
 
 The map $\Lam : \nu_*^{\gamma,\si}\in\{\trm{$*$-infinitely divisible laws}\} \to  \nu_\bxp^{\gamma,\si} \in \{\trm{$\bxp$-infinitely divisible laws}\}$ is called the {\it Bercovici-Pata bijection}. In \cite{appenice} (see also \cite{gotzea}), Bercovici and Pata proved that this bijection preserves limit theorems, i.e. that for any sequence $(k_n)$ of positive integers tending to infinity, for any sequence $(\nu_n)$ of laws on $\R$, for any $*$-infinitely divisible law $\mu$, we have $$  \underbrace{\nu_n*\cdots*\nu_n}_{\trm{$k_n$ times}}\lto \mu\iff 
 \underbrace{\nu_n\bxp\cdots\bxp\nu_n}_{\trm{$k_n$ times}}\lto \La(\mu).$$

 \subsubsection{The cumulants point of view}\la{introcum}
Let $\mu$ be a \pro measure on $\R$   whose Laplace  transform is defined in a neighborhood of zero. Its {\it classical cumulants} are the numbers  $(c_n(\mu))_{n\ge 1}$ defined by the formula $$ \qquad\qquad  \qquad\qquad \log \int_{t\in \R} e^{\xi t}\ud \mu(t) =\sum_{n\ge 1}\f{c_n(\mu)}{n!}\xi^n \qquad\qquad \trm{($\xi\in \C$ small enough).}$$

In the same way, for $\mu$ a compactly supported \pro measure on $\R$, the 
{\it free cumulants} of $\mu$ are the numbers $(k_n(\mu))_{n\ge 1}$ defined by the formula $$ \qquad\qquad  \qquad\qquad R_\mu(z) =\sum_{n\ge 1}k_n(\mu)z^{n-1} \qquad\qquad \trm{($z\in \C$ small enough).}$$

It can easily be seen (see e.g. \cite[Eq. (2.1)]{FloTaylor}) that for any \Lvy pair $(\gamma, \mu)$, the classical (resp. free) cumulants of $\nu_*^{\gamma, \si}$ (resp. of $\nu_\bxp^{\gamma, \si}$) are the given by the formula \be\la{571119h57}k_n(\nu_*^{\gamma, \si})=c_n(\nu_\bxp^{\gamma, \si})=\beg{cases}\gamma+\int t\ud \si(t)&\trm{ if $n=1$,}\\
\\
\int t^{n-2}(1+t^2)\ud \si(t)&\trm{ if $n\ge 2$.}\enc\ee
Hence the Bercovici-Pata bijection can be seen as (the continuous extension of) the  map which transforms classical cumulants into free ones. 

\subsection{Combinatorial definitions}\label{appendicecombi}




 
In this section, definitions and results can be found in the classical book \cite{MR0357171} of Berge.

\subsubsection{Graphs}\la{section_graphs}
A \emph{graph} $G=(V,E)$ is characterized by a set $V$ of vertices and a family $E=(e_i, i\in I)$ of edges:
\begin{enumerate}
\item If the edges are elements of $V\times V$, $G$ is a \emph{directed} graph;
\item If the edges are singletons or pairs of $V$, $G$ is a \emph{non-directed} graph.
\end{enumerate}
Note that {\bf multiple edges} are possible with this definition, since 
$E=(e_i, i\in I)$ is a \emph{family} and not a {\it set}, which means that the $e_i$'s are not pairwise distinct.
\\

In the following three definitions, we suppose {\bf $G$  non-directed}.
\begin{enumerate}
\item A \emph{connected component} of $G$ is a block of the smallest partition of $V$ which is above $E$ w.r.t. the refinement order.
\item A {\it path} in $G$ is a sequence  $(v_0,e_{i_1},v_1,e_{i_2},\ldots,e_{i_k},v_{k})$ with $k\geq 1$,  $v_0,\ldots,v_{k}\in V$, $e_{i_j}=\{v_{j-1},v_{j}\}$. If, moreover, the $v_0=v_k$, then the path is said to be a {\it closed}. In this case, if $v_0, \ld, v_{k-1}$ are pairwise distinct and $i_1, \ld, i_k$  are pairwise distinct, the path is said to be a {\it cycle}. 
\item The \emph{cyclomatic number} of $G$ is equal to 
$
\vert I\vert-\vert V\vert+n
$,
with $n$ the number of connected components of $G$.\\
\end{enumerate}

The following proposition can easily be proved by induction on the number
of vertices of $G$. 

\begin{propo}\label{ptn:0513}
For every non directed  graph $G$, its cyclomatic number is non-negative. It vanishes if and only if $G$ has no cycle.
\end{propo}

The previous definitions extend easily to the framework of a {\bf directed} graph
 $G=(V,E=(e_i, i\in I))$. For example, a {\it path} in $G$ is a sequence  $(v_0,e_{i_1},v_1,e_{i_2},\ldots,e_{i_k},v_{k})$ with $\ell\geq 1$,  $v_0,\ldots,v_{k}\in V$, $e_{i_j}=(v_{j-1},v_{j})$. 
In a directed graph,  two cycles $c,c'$  are said to be {\it disjoint} if they do not have any edge in common, \ie if there is no $i\in I$ \st $e_i$ appears both in $c$ and in $c'$ (whereas they can have some vertices in common). 
We say that the directed graph $G$ is a {\it circuit} if a closed path of $G$ visits each of its edges exactly once (vertices can be passed-by more than once). One can easily prove (by induction) that any circuit is a union of  disjoint  cycles.\\

\subsubsection{Hypergraphs}
A {\it hypergraph} is  a pair  $H=(V,E)$ , where   $V$  is a set (the \emph{vertices}) and
$E=(E_i,i\in I)$ is a family (the \emph{edges}), of non-empty subsets of $V$ such
that $\cup_{i\in I}E_i=V$. 
An example is given at Figure \re{ex_hypergraph}. 
 If the edges have
only one or two elements, then $H$ reduces to a non directed graph. 
 
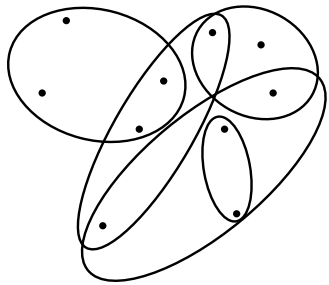
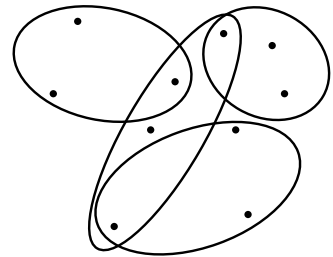
\begin{figure}[h!]\centering
\subfigure[A hypergraph]{
\scalebox{0.8} 
{
\begin{pspicture}(0,-2.4322512)(5.759434,2.5269399)
\psdots[dotsize=0.12](3.6335897,1.9466037)
\psdots[dotsize=0.12](4.4335895,1.7466037)
\psdots[dotsize=0.12](4.6335897,0.94660366)
\psdots[dotsize=0.12](3.8335898,0.3466037)
\psdots[dotsize=0.12](2.4335897,0.3466037)
\psdots[dotsize=0.12](1.2335896,2.1466036)
\psdots[dotsize=0.12](0.8335897,0.94660366)
\psdots[dotsize=0.12](1.8335897,-1.2533963)
\psdots[dotsize=0.12](4.03359,-1.0533963)
\psdots[dotsize=0.12](2.8335898,1.1466037)
\rput{-30.0}(-0.14271091,2.360603){\psellipse[linewidth=0.04,dimen=outer](4.3335896,1.4466037)(1.1,0.9)}
\rput{-30.0}(0.20369925,1.3756665){\psellipse[linewidth=0.04,dimen=outer](2.6688783,0.30772528)(0.66961527,2.2428203)}
\rput{40.0}(0.55725,-2.3393438){\psellipse[linewidth=0.04,dimen=outer](3.4922721,-0.40415594)(2.5,0.9985081)}
\rput{-15.0}(-0.26357415,0.49116302){\psellipse[linewidth=0.04,dimen=outer](1.7335896,1.2466037)(1.5,1.1)}
\rput{10.0}(0.004427832,-0.677403){\psellipse[linewidth=0.04,dimen=outer](3.8735898,-0.3133963)(0.4,0.9)}
\end{pspicture} 
}\la{ex_hypergraph}
}
\hspace{3cm}
\subfigure[A hypergraph with no cycle]{
\scalebox{0.8}
{
\begin{pspicture}(0,-2.1950872)(5.552422,2.2299614)
\psdots[dotsize=0.12](3.6541095,1.6398984)
\psdots[dotsize=0.12](4.454109,1.4398984)
\psdots[dotsize=0.12](4.6541095,0.6398984)
\psdots[dotsize=0.12](3.8541095,0.03989839)
\psdots[dotsize=0.12](2.4541094,0.03989839)
\psdots[dotsize=0.12](1.2541094,1.8398983)
\psdots[dotsize=0.12](0.8541094,0.6398984)
\psdots[dotsize=0.12](1.8541094,-1.5601016)
\psdots[dotsize=0.12](4.0541096,-1.3601016)
\psdots[dotsize=0.12](2.8541095,0.8398984)
\rput{-30.0}(0.013390855,2.3297722){\psellipse[linewidth=0.04,dimen=outer](4.3541093,1.1398984)(1.1,0.9)}
\rput{-30.0}(0.35980102,1.3448356){\psellipse[linewidth=0.04,dimen=outer](2.689398,0.0010199854)(0.66961527,2.2428203)}
\rput{20.0}(-0.121389285,-1.1600785){\psellipse[linewidth=0.04,dimen=outer](3.2288713,-0.92425567)(1.7740905,0.9985083)}
\rput{-12.0}(-0.19926068,0.37075305){\psellipse[linewidth=0.04,dimen=outer](1.6641095,1.1332959)(1.5,0.94697905)}
\end{pspicture} 
}\label{fig:hypergraph_with_no_cycle}
}\caption{Two examples of  hypergraphs}
\end{figure}

The refinement order for hypergraphs defined on a same set of vertices is defined as the refinement order on their sets of edges.

A {\it connected component} of $H$ is a block of the smallest
  partition of $V$ which is above $E$ w.r.t. the refinement order. For
  instance, there is only one connected component in the 
  hypergraph of Figure \re{ex_hypergraph}, hence it is said to be  {\it  connected}.

Cycles are defined in hypergraphs as in non-directed  graphs, except that cycles with length one are not accepted (in fact, they do not make sense) : 
a {\it cycle} of $H$ is a sequence $(v_0,E_{i_1},v_1,E_{i_2},\ldots,v_{\ell-1},E_{i_\ell},v_{\ell})$ with $\ell\geq 2$,  $v_0,v_1,\cdots,v_{\ell-1}$ pairwise distincts, $v_0=v_{\ell}$, $i_1,\ldots,i_\ell$ pairwise distinct and $v_{j-1},v_{j}\in E_{i_j}$ for all $j$.

An example of hypergraph with no cycle is given in Figure \re{fig:hypergraph_with_no_cycle}.
Notice that a hypergraph with no cycle is linear: two different edges have at most one vertex in common. Notice also that a hypergraph with only one vertex or only one edge has no cycle.
\begin{propo}\label{prop:0503}
Let $n$ be the number of connected components of $H$. Then:

$\bullet$ $\sum_{i\in I}\vert E_i\vert-\vert I\vert-\vert V\vert+n\geq 0$;

$\bullet$ $\sum_{i\in I}\vert E_i\vert-\vert I\vert-\vert V\vert+n= 0$ if and only if $H$ has no cycle.
\end{propo}

\beg{pr}This is a standard result (see for instance \cite[Ch. 7, Prop. 4]{MR0357171}), which is a consequence of Proposition \ref{ptn:0513}. For the easyness of the reader, let us recall the short proof. We associate to $H$ a 
non directed graph $G(H)$ with vertices $V\cup I$, and with an edge between $x\in V$ and $i\in I$ if and only if $x\in E_i$. The graph  $G(H)$ is simple, i.e. there cannot be more than one edge between two of its vertices. Also, it has no loop, i.e. no edge from one vertex to itself. The number of edges is equal to $ \sum_{i\in I}\vert E_i\vert$, and the number of connected components is $n$, so the cyclomatic number of $G(H)$ is equal to
$
\sum_{i\in I}\vert E_i\vert-(\vert V\vert+\vert I\vert)+n
$.
It is non-negative, and vanishes if and only $G(H)$ has no cycle, which means exactly that $H$ has no cycle.\en{pr}


\end{document}